\documentstyle[titlepage,12pt]{article}

\newcommand{\Ha}{\put(0,30){\circle*{3}}}
\newcommand{\Hb}{\put(-30,15){\circle*{3}}}
\newcommand{\Hc}{\put(-30,-15){\circle*{3}}}
\newcommand{\Hd}{\put(0,-30){\circle*{3}}}
\newcommand{\He}{\put(30,-15){\circle*{3}}}

\newcommand{\Hab}{\put(0,30){\line(-2,-1){30}}}

\newcommand{\Hae}{\put(0,30){\line(2,-3){30}}}

\newcommand{\Hbc}{\put(-30,15){\line(0,-1){30}}}

\newcommand{\Hcd}{\put(-30,-15){\line(2,-1){30}}}

\newcommand{\Hde}{\put(0,-30){\line(2,1){30}}}

\newcommand{\Gca}{\put(20,0){\circle*{3}}}
\newcommand{\Gda}{\put(30,0){\circle*{3}}}
\newcommand{\Gea}{\put(40,0){\circle*{3}}}

\newcommand{\Gac}{\put(0,20){\circle*{3}}}

\newcommand{\Gdc}{\put(30,20){\circle*{3}}}

\newcommand{\Ggc}{\put(60,20){\circle*{3}}}

\newcommand{\Gae}{\put(0,40){\circle*{3}}}

\newcommand{\Gge}{\put(60,40){\circle*{3}}}

\newcommand{\Gcg}{\put(20,60){\circle*{3}}}
\newcommand{\Gdg}{\put(30,60){\circle*{3}}}
\newcommand{\Geg}{\put(40,60){\circle*{3}}}

\newcommand{\GaaV}[1]{\put(0,0){\circle*{#1}}}
\newcommand{\GbaV}[1]{\put(10,0){\circle*{#1}}}

\newcommand{\GdaV}[1]{\put(30,0){\circle*{#1}}}
\newcommand{\GeaV}[1]{\put(40,0){\circle*{#1}}}

\newcommand{\GabV}[1]{\put(0,10){\circle*{#1}}}

\newcommand{\GcbV}[1]{\put(20,10){\circle*{#1}}}
\newcommand{\GdbV}[1]{\put(30,10){\circle*{#1}}}

\newcommand{\GfbV}[1]{\put(50,10){\circle*{#1}}}
\newcommand{\GgbV}[1]{\put(60,10){\circle*{#1}}}

\newcommand{\GacV}[1]{\put(0,20){\circle*{#1}}}
\newcommand{\GbcV}[1]{\put(10,20){\circle*{#1}}}
\newcommand{\GccV}[1]{\put(20,20){\circle*{#1}}}

\newcommand{\GecV}[1]{\put(40,20){\circle*{#1}}}

\newcommand{\GgcV}[1]{\put(60,20){\circle*{#1}}}

\newcommand{\GaeV}[1]{\put(0,40){\circle*{#1}}}

\newcommand{\GceV}[1]{\put(20,40){\circle*{#1}}}
\newcommand{\GdeV}[1]{\put(30,40){\circle*{#1}}}
\newcommand{\GeeV}[1]{\put(40,40){\circle*{#1}}}
\newcommand{\GfeV}[1]{\put(50,40){\circle*{#1}}}
\newcommand{\GgeV}[1]{\put(60,40){\circle*{#1}}}

\newcommand{\GdgV}[1]{\put(30,60){\circle*{#1}}}

\newcommand{\GcaL}[2]{\Gca \put(20,0){\makebox(0,0){\hspace{#1}#2}}}

\newcommand{\GeaL}[2]{\Gea \put(40,0){\makebox(0,0){\hspace{#1}#2}}}

\newcommand{\GacL}[2]{\Gac \put(0,20){\makebox(0,0){\hspace{#1}#2}}}

\newcommand{\GdcL}[2]{\Gdc \put(30,20){\makebox(0,0){\hspace{#1}#2}}}

\newcommand{\GgcL}[2]{\Ggc \put(60,20){\makebox(0,0){\hspace{#1}#2}}}

\newcommand{\GaeL}[2]{\Gae \put(0,40){\makebox(0,0){\hspace{#1}#2}}}

\newcommand{\GgeL}[2]{\Gge \put(60,40){\makebox(0,0){\hspace{#1}#2}}}

\newcommand{\GcgL}[2]{\Gcg \put(20,60){\makebox(0,0){\hspace{#1}#2}}}

\newcommand{\GegL}[2]{\Geg \put(40,60){\makebox(0,0){\hspace{#1}#2}}}

\newcommand{\GaaLV}[3]{\GaaV{#3} \put(0,0){\makebox(0,0){\hspace{#1}#2}}}

\newcommand{\GacLV}[3]{\GacV{#3} \put(0,20){\makebox(0,0){\hspace{#1}#2}}}
\newcommand{\GbcLV}[3]{\GbcV{#3} \put(10,20){\makebox(0,0){\hspace{#1}#2}}}
\newcommand{\GccLV}[3]{\GccV{#3} \put(20,20){\makebox(0,0){\hspace{#1}#2}}}

\newcommand{\GecLV}[3]{\GecV{#3} \put(40,20){\makebox(0,0){\hspace{#1}#2}}}

\newcommand{\GgcLV}[3]{\GgcV{#3} \put(60,20){\makebox(0,0){\hspace{#1}#2}}}

\newcommand{\GaeLV}[3]{\GaeV{#3} \put(0,40){\makebox(0,0){\hspace{#1}#2}}}

\newcommand{\GceLV}[3]{\GceV{#3} \put(20,40){\makebox(0,0){\hspace{#1}#2}}}
\newcommand{\GdeLV}[3]{\GdeV{#3} \put(30,40){\makebox(0,0){\hspace{#1}#2}}}
\newcommand{\GeeLV}[3]{\GeeV{#3} \put(40,40){\makebox(0,0){\hspace{#1}#2}}}
\newcommand{\GfeLV}[3]{\GfeV{#3} \put(50,40){\makebox(0,0){\hspace{#1}#2}}}
\newcommand{\GgeLV}[3]{\GgeV{#3} \put(60,40){\makebox(0,0){\hspace{#1}#2}}}

\newcommand{\Gaaac}{\put(0,0){\line(0,1){20}}}

\newcommand{\Gaacc}{\put(0,0){\line(1,1){20}}}

\newcommand{\Gbaab}{\put(10,0){\line(-1,1){10}}}

\newcommand{\Gbacb}{\put(10,0){\line(1,1){10}}}

\newcommand{\Gcaea}{\put(20,0){\line(1,0){20}}}

\newcommand{\Gcaac}{\put(20,0){\line(-1,1){20}}}

\newcommand{\Gcaeg}{\put(20,0){\line(1,3){20}}}

\newcommand{\Gdaac}{\put(30,0){\line(-3,2){30}}}
\newcommand{\Gdabc}{\put(30,0){\line(-1,1){20}}}
\newcommand{\Gdacc}{\put(30,0){\line(-1,2){10}}}
\newcommand{\Gdadc}{\put(30,0){\line(0,1){20}}}
\newcommand{\Gdaec}{\put(30,0){\line(1,2){10}}}

\newcommand{\Gdagc}{\put(30,0){\line(3,2){30}}}

\newcommand{\Gdaae}{\put(30,0){\line(-3,4){30}}}

\newcommand{\Gdage}{\put(30,0){\line(3,4){30}}}

\newcommand{\Geadb}{\put(40,0){\line(-1,1){10}}}

\newcommand{\Geafb}{\put(40,0){\line(1,1){10}}}

\newcommand{\Geacg}{\put(40,0){\line(-1,3){20}}}

\newcommand{\Geaeg}{\put(40,0){\line(0,1){60}}}

\newcommand{\Gacae}{\put(0,20){\line(0,1){20}}}

\newcommand{\Gacce}{\put(0,20){\line(1,1){20}}}

\newcommand{\Gacee}{\put(0,20){\line(2,1){40}}}
\newcommand{\Gacfe}{\put(0,20){\line(5,2){50}}}

\newcommand{\Gaccg}{\put(0,20){\line(1,2){20}}}

\newcommand{\Gbcae}{\put(10,20){\line(-1,2){10}}}

\newcommand{\Gbcde}{\put(10,20){\line(1,1){20}}}

\newcommand{\Gbcge}{\put(10,20){\line(5,2){50}}}

\newcommand{\Gccae}{\put(20,20){\line(-1,1){20}}}

\newcommand{\Gccce}{\put(20,20){\line(0,1){20}}}
\newcommand{\Gccde}{\put(20,20){\line(1,2){10}}}
\newcommand{\Gccee}{\put(20,20){\line(1,1){20}}}

\newcommand{\Gdcae}{\put(30,20){\line(-3,2){30}}}

\newcommand{\Gdcge}{\put(30,20){\line(3,2){30}}}

\newcommand{\Gecde}{\put(40,20){\line(-1,2){10}}}

\newcommand{\Gecfe}{\put(40,20){\line(1,2){10}}}

\newcommand{\Ggcfe}{\put(60,20){\line(-1,2){10}}}
\newcommand{\Ggcge}{\put(60,20){\line(0,1){20}}}

\newcommand{\Ggccg}{\put(60,20){\line(-1,1){40}}}

\newcommand{\Ggceg}{\put(60,20){\line(-1,2){20}}}

\newcommand{\Gaedg}{\put(0,40){\line(3,2){30}}}

\newcommand{\Gcedg}{\put(20,40){\line(1,2){10}}}

\newcommand{\Gdedg}{\put(30,40){\line(0,1){20}}}

\newcommand{\Gfedg}{\put(50,40){\line(-1,1){20}}}

\newcommand{\Ggecg}{\put(60,40){\line(-2,1){40}}}
\newcommand{\Ggedg}{\put(60,40){\line(-3,2){30}}}
\newcommand{\Ggeeg}{\put(60,40){\line(-1,1){20}}}

\newcommand{\ben}{\begin{enumerate}}
\newcommand{\een}{\end{enumerate}}
\newcommand{\ble}{\begin{lem}}
\newcommand{\ele}{\end{lem}}
\newcommand{\bth}{\begin{thm}}
\newcommand{\eth}{\end{thm}}
\newcommand{\bpr}{\begin{prop}}
\newcommand{\epr}{\end{prop}}
\newcommand{\bco}{\begin{cor}}
\newcommand{\eco}{\end{cor}}
\newcommand{\bcon}{\begin{conj}}
\newcommand{\econ}{\end{conj}}
\newcommand{\bde}{\begin{defn}}
\newcommand{\ede}{\end{defn}}
\newcommand{\bex}{\begin{exa}}
\newcommand{\eex}{\end{exa}}
\newcommand{\barr}{\begin{array}}
\newcommand{\earr}{\end{array}}
\newcommand{\btab}{\begin{tabular}}
\newcommand{\etab}{\end{tabular}}
\newcommand{\beq}{\begin{equation}}
\newcommand{\eeq}{\end{equation}}
\newcommand{\bea}{\begin{eqnarray*}}
\newcommand{\eea}{\end{eqnarray*}}
\newcommand{\bce}{\begin{center}}
\newcommand{\ece}{\end{center}}
\newcommand{\bpi}{\begin{picture}}
\newcommand{\epi}{\end{picture}}
\newcommand{\bfi}{\begin{figure} \begin{center}}
\newcommand{\efi}{\end{center} \end{figure}}
\newcommand{\capt}{\caption}
\newcommand{\bsl}{\begin{slide}{}}
\newcommand{\esl}{\end{slide}}

\newcommand{\bib}{thebibliography}
\newcommand{\pf}{{\bf Proof.}}
\newcommand{\qed}{\rule{1ex}{1ex}}
\newcommand{\Qed}{\rule{1ex}{1ex} \medskip}

\newcommand{\rp}[2]{\rule{#1pt}{#2pt}}

\newcommand{\mbl}[1]{\makebox(0,0)[l]{#1}}
\newcommand{\mbr}[1]{\makebox(0,0)[r]{#1}}
\newcommand{\mbt}[1]{\makebox(0,0)[t]{#1}}
\newcommand{\mbb}[1]{\makebox(0,0)[b]{#1}}

\newcommand{\emp}{\emptyset}
\newcommand{\sbs}{\subset}
\newcommand{\sbe}{\subseteq}

\newcommand{\spe}{\supseteq}
\newcommand{\setm}{\setminus}

\newcommand{\iso}{\cong}
\newcommand{\ptn}{\vdash}

\newcommand{\zh}{\hat{0}}
\newcommand{\oh}{\hat{1}}

\newcommand{\lt}{\lhd}

\newcommand{\lte}{\unlhd}
\newcommand{\gte}{\unrhd}
\newcommand{\jn}{\vee}
\newcommand{\Jn}{\bigvee}
\newcommand{\mt}{\wedge}

\newcommand{\sd}{\bigtriangleup}
\newcommand{\case}[4]{\left\{\barr{ll}#1&\mbox{#2}\\#3&\mbox{#4}\earr\right.}

\newcommand{\ra}{\rightarrow}

\newcommand{\al}{\alpha}
\newcommand{\be}{\beta}
\newcommand{\ga}{\gamma}
\newcommand{\de}{\delta}
\newcommand{\ep}{\epsilon}
\newcommand{\io}{\iota}
\newcommand{\la}{\lambda}
\newcommand{\mut}{\tilde{\mu}}

\newcommand{\si}{\sigma}

\newcommand{\De}{\Delta}

\newcommand{\ba}{{\bf a}}
\newcommand{\bb}{{\bf b}}

\newcommand{\bs}{{\bf s}}

\newcommand{\bu}{{\bf u}}
\newcommand{\bv}{{\bf v}}
\newcommand{\bw}{{\bf w}}
\newcommand{\bx}{{\bf x}}
\newcommand{\by}{{\bf y}}
\newcommand{\bA}{{\bf A}}

\newcommand{\bZ}{{\bf Z}}

\newcommand{\cA}{{\cal A}}
\newcommand{\cB}{{\cal B}}

\newcommand{\cL}{{\cal L}}

\newcommand{\cP}{{\cal P}}

\newcommand{\cS}{{\cal S}}

\newcommand{\cW}{{\cal W}}

\newcommand{\dm}{Discrete Math.}

\newcommand{\ejc}{European J. Combin.}

\newcommand{\jams}{J. Amer. Math. Soc.}

\newcommand{\jcta}{J. Combin. Theory Ser. A}

\newtheorem{thm}{Theorem}[section]
\newtheorem{prop}[thm]{Proposition}
\newtheorem{cor}[thm]{Corollary}
\newtheorem{lem}[thm]{Lemma}
\newtheorem{conj}[thm]{Conjecture}
\newtheorem{exa}[thm]{Example}
\newtheorem{defn}[thm]{Definition}

\begin{document}
\pagestyle{empty}
\title{M\"obius functions of lattices}
\author{Andreas Blass\\
Department of Mathematics \\ University of Michigan \\
Ann Arbor, MI 48109-1003\\
ablass@umich.edu\\[5pt]
and\\[5pt]
Bruce E. Sagan \\
Department of Mathematics \\ Michigan State University
\\ East Lansing, MI 48824-1027\\sagan@math.msu.edu}

\date{\today \\[1in]
	\begin{flushleft}
	\end{flushleft}
       }
\maketitle

\begin{flushleft} Proposed running head: \end{flushleft}
	\begin{center} 
M\"obius functions of lattices
	\end{center}

Send proofs to:
\begin{center}
Bruce E. Sagan \\ Department of Mathematics \\Michigan State
University \\ East Lansing, MI 48824-1027\\[5pt]
Tel.: 517-355-8329\\
FAX: 517-336-1562\\
Email: sagan@mth.msu.edu
\end{center}

	\begin{abstract} We introduce the concept of a bounded below
set in a lattice.  This can be used to give a generalization of Rota's
broken circuit theorem to any finite lattice.  We then show how this
result can be used to compute and combinatorially explain the M\"obius
function in various examples including non-crossing set partitions,
shuffle posets, and integer partitions in dominance order.  Next we present a 
generalization of Stanley's theorem that the characteristic polynomial
of a semimodular supersolvable lattice factors over the integers.  We
also give some applications of this second main theorem, including the
Tamari lattices.\end{abstract}
\pagestyle{plain}

\section{Bounded below sets}			\label{bbs}

In a fundamental paper~\cite{whi:lem}, Whitney showed how broken circuits
could be used to compute the coefficients of the chromatic polynomial
of a graph.  
In another seminal paper~\cite{rot:tmf}, Rota
refined and extended Whitney's theorem to give a characterization of
the M\"obius function of a geometric lattice.  Then one of
us~\cite{sag:grn} generalized Rota's result to a larger class of
lattices.  In this paper we will present a theorem for an arbitrary
finite lattice
that includes all the others as special cases.  To do so, we shall
need 
to replace the notion of a broken circuit by a new one which we call a
bounded below set.  Next we present some applications to lattices
whose M\"obius functions had previously been computed but in a less
simple or less combinatorial way:  shuffle posets~\cite{gre:ps},
non-crossing set partition lattices~\cite{kre:pnc,rei:npc}, and
integer partitions under dominance
order~\cite{bog:mfd,bry:lip,gre:clm}.  

The second half of the paper is dedicated to generalizing the result
of Stanley~\cite{sta:ssl} that the characteristic polynomial of a
semimodular supersolvable lattice factors over the integers.
We replace both supersolvability and semimodularity by weaker
conditions which we call left-modularity and the level condition,
respectively, in such a way that the conclusion still holds.
Examples of this factorization (not covered by Stanley's theorem) are
provided by certain shuffle posets 
and the Tamari lattices~\cite{ft:tfi,ht:spl}.
We end with a further generalization of Rota's theorem and a section
of comments and questions. 

Throughout this paper $L$ will denote a finite lattice.  Any relevant
definitions not given can be found in Stanley's text~\cite{sta:ec1}.
We will use $\mt$ for the meet (greatest lower bound) and $\jn$ for
the join (least upper bound) in $L$.  Since $L$ is finite it also has
a unique minimal element $\zh$ and a unique maximal element $\oh$.
The {\sl M\"obius function of $L$}, $\mu: L\ra \bZ$, is defined
recursively by
	$$\mu(x)=\case{1}{if $x=\zh$,}{-\sum_{y<x}\mu(y)}{if $x>\zh$.}$$
We let $\mu(L)=\mu(\oh)$.  Note that $\mu$ is the unique $\bZ$-valued
function on $L$ such that $\sum_{y\le x}\mu(y)=\de_{\zh x}$ (Kronecker
delta).  

Our goal is to give a new combinatorial description of $\mu(x)$.  Let
$A(L)$ be the set of atoms of $L$, i.e., those elements covering
$\zh$.
Give $A(L)$ an arbitrary partial order, which we denote $\lte$ to
distinguish it from the partial order $\le$ in $L$.  So $\lte$ can be
anything from a total order to the total incomparability order induced
by $\le$.  A nonempty set $D\sbe A(L)$ is {\sl bounded below} or {\sl
BB} if, for every 
$d\in D$ there is an $a\in A(L)$ such that
	\begin{eqnarray}
	a&\lt& d \mbox{ and }		\label{BB1}\\
	a&<&\Jn D.			\label{BB2}
	\end{eqnarray}
So $a$ is simultaneously a strict lower bound for $d$ in the order $\lte$ and 
for $\Jn D$ in $\le$.  We will say that $B\sbe A(L)$ is {\sl NBB} 
if $B$ does not contain any $D$ which is bounded
below.  In this case we will call $B$ an {\sl NBB base} for $x=\Jn B$.
We can now state our main result.
\bth							\label{NBB}
Let $L$ be any finite lattice and let $\lte$ be any partial order on $A(L)$.
Then for all $x\in L$ we have
	\beq						\label{NBBsum}
	\mu(x)=\sum_B (-1)^{|B|}
	\eeq
where the sum is over all NBB bases $B$ of $x$ and $|\cdot|$ denotes
cardinality.
\eth
\pf\  For $x\in L$ define $\mut(x)=\sum_B (-1)^{|B|}$ summed over all
NBB bases $B$ of $x$.  To prove that $\mut(x)=\mu(x)$ it suffices to
show $\sum_{y\le x}\mut(y)=\de_{\zh x}$.  If $x=\zh$ then $x=\Jn B$
only for $B=\emp$ which is NBB.  So
	$$\sum_{y\le \zh}\mut(y)=\mut(\zh)=(-1)^0=1$$
as desired.

If $x>\zh$ then to get $\sum_{y\le x}\mut(y)=0$ we first set up a
corresponding signed set $\cS$.
Let
	$$\cS=\{B\ :\ B \mbox{ an NBB base for some } y\le x\}.$$	
The sign of $B\in\cS$ will be $\ep(B)=(-1)^{|B|}$.  Then from the definitions
	$$\sum_{y\le x}\mut(y)=\sum_{B\in\cS} \ep(B).$$
If we can find a sign-reversing involution on $\cS$ then this last sum
will be zero and we will be done.
{}From the set of atoms $a\le x$ pick one, $a_0$,
which is minimal with respect to $\lte$.   Consider the map on $\cS$
defined by $\io(B)=B\sd a_0$ where $\sd$ is symmetric difference.
(Here and below we omit 
set braces around singletons writing, for example, $a_0$ instead of 
$\{a_0\}$.) 
Clearly $\io$ is a sign-reversing involution provided it is well
defined, i.e., we must check that $B$ being NBB implies that  $\io(B)$
is NBB.

If $\io(B)=B\setm a_0$ then clearly $\io(B)$ is still NBB.  We
will do the case $B':=\io(B)=B\cup a_0$ by contradiction.  Suppose
$B'\spe D$ where $D$ is bounded below.  Then $a_0\in D$ since $B$
itself is NBB.  Let $a$ be the corresponding element guaranteed by the
definition of a bounded below set. Then $a\lt a_0$ and $a<\Jn B'\le x$
which contradicts the definition of $a_0$.\hfill\qed

\thicklines
\setlength{\unitlength}{1.3pt}
\bfi
\btab{cc}
\bpi(90,80)(-30,-10)
\put(-50,30){$(L,\le)=$}
\Gda
\put(30,-10){\makebox(0,0){$\zh$}}
\GacL{-15pt}{$a$}
\GdcL{-15pt}{$b$}
\GgcL{15pt}{$c$}
\GaeL{-15pt}{$x$}
\GgeL{15pt}{$y$}
\Gdg
\put(30,70){\makebox(0,0){$\oh$}}
\Gdaac \Gdadc \Gdagc \Gacae \Gdcae \Gdcge \Ggcge \Gaedg \Ggedg
\epi
&
\bpi(150,80)(-90,-10)
\put(-70,30){$(A(L),\lte)=$}
\Gda
\put(30,-10){\makebox(0,0){$b$}}
\GaeL{-15pt}{$a$}
\GgeL{15pt}{$c$}
\Gdaae \Gdage
\epi
\etab
\capt{A lattice $L$ and partial order $\lte$ on $A(L)$}	\label{l}	
\efi

Here is  an example to illustrate this result.  Suppose the lattice
$L$ and partial order $\lte$ are as given in Figure~\ref{l}.  To find
the bounded below sets, note that by~(\ref{BB1}) no set containing an
element minimal in $\lte$ is BB.  Furthermore~(\ref{BB1})
and~(\ref{BB2}) together imply that no single element set is BB
either.  (These observations will be important in our other examples.)
Thus the only possible BB set is $\{a,c\}$, and it does satisfy the
definition since $b\lt a,c$ and $b<\Jn\{a,c\}=\oh$.
So $x$ has one NBB base, namely $\{a,b\}$ and so $\mu(x)=(-1)^2$ which
is easily checked from the definition of $\mu$.  Similarly $\mu(y)=(-1)^2$.
Finally $\oh$ has no NBB bases and so $\mu(\oh)=0$ (the empty sum).

We should see why our theorem implies Rota's broken circuit result.
To do this we need to recall some definitions.
Let $L$ be a geometric lattice with rank function $\rho$.  It is
well known, and easy to prove, that if $B\sbe A(L)$ then 
$\rho(\Jn B)\le|B|$.  Call $B$ {\sl independent} if $\rho(\Jn B)=|B|$
and {\sl dependent} otherwise.  A {\sl circuit} $C$ is a minimal (with
respect to inclusion) dependent set.  Now let $\lte$ be any total order
on $A(L)$.  Then each circuit $C$ gives rise to a {\sl broken circuit}
$C'=C\setm c$ where $c$ is the first element of $C$ under $\lte$.  A
set $B\sbe A(L)$ is {\sl NBC} (no broken circuit) if $B$ does not
contain any broken circuit and in this case $B$ is an {\sl NBC base}
for $x=\Jn B$.  Rota's NBC theorem~\cite{rot:tmf} is as follows.
\bth[Rota]					\label{NBC}
Let $L$ be a finite geometric lattice and let $\lte$ be any total
order on $A(L)$. 
Then for all $x\in L$ we have
    $$\mu=(-1)^{\rho(x)}\cdot(\mbox{number of NBC bases of $x$).}$$
\eth

To derive this result from Theorem~\ref{NBB}, we first prove that when $L$ is
geometric and $\lte$ is total then the NBB and NBC sets coincide.  For
this it suffices to show that every broken circuit is bounded below
and that every bounded below set contains a broken circuit.  If
$C'=C\setm c$ is a broken circuit, then for every $c'\in C'$ we have
$c\lt c'$ and $c<\Jn C=\Jn C'$ so $C'$ is BB.  For the other 
direction, if $D$ is bounded below then consider the $\lte$-first
element $d$ of $D$ and let $a\in A(L)$ be
the element guaranteed by the BB definition.  Then by~(\ref{BB2}) we have
$\rho(\Jn D\jn a)<|D\cup a|$.  So $D\cup a$ is dependent and
contains a circuit $C$.  Now~(\ref{BB1}) and the choice of $a$ and $d$ show
that for the corresponding 
broken circuit we have $C'\sbe D$.  Thus NBB and NBC sets are the same
in this setting.
Finally Rota's expression for $\mu(x)$ is obtained from ours by noting
that when $L$ is semimodular then all NBB bases for $x$ have the same
size, namely $\rho(x)$.  Similar arguments show that the main result
of~\cite{sag:grn} is a special case of Theorem~\ref{NBB}.

Another corollary of this theorem is a special case of Rota's Crosscut
Theorem~\cite{rot:tmf}.
\bth[Rota]						\label{cct}	
Let $L$ be any finite lattice.  Then for all $x\in L$
	\beq							\label{cross}
	\mu(x)=\sum_B (-1)^{|B|}
	\eeq
where the sum is over all $B\sbe A(L)$ such that $\Jn B=x$.
\eth
To obtain~(\ref{cross}) it suffices to take $\lte$ to be the total
incomparability order in Theorem~\ref{NBB}.  In fact our theorem (where
the partial order on $A(L)$ is arbitrary) can
be viewed as interpolating between the Crosscut Theorem (where 
$A(L)$ forms an antichain) and the NBC Theorem (where $A(L)$ forms a chain).

This raises the question of why one would want to consider an
arbitrary partial order $\lte$ on $A(L)$ when one can always take the one
induced by the order in $L$.  The reason is that the number of terms
in the sum~(\ref{cross}) is generally much larger than the number
in~(\ref{NBBsum}).  From the viewpoint of  efficient computation
of $\mu$, the best scenario is the same as the one in the geometric
case where~(\ref{NBBsum}) has exactly $|\mu(x)|$ terms, all of the
same sign.  A partial order $\lte$ on $A(L)$ for which
this happens for all $x\in L$ will be  called {\sl perfect}.  There
are posets where no such $\lte$ exists, such as the $k$-equal
intersection lattices which were introduced by Bj\"orner, Lov\'asz and
Yao~\cite{bl:ldt,bly:ldt}.  However, all the examples we will consider
in the next sections are perfect.  Another thing to note is that if
$\lte$ is  perfect then so is any linear extension of it.  However, to
make the combinatorics of $\mu$ as clear as
possible it is often best to take a perfect $\lte$ with the 
least possible number of order relations.

\section{Non-crossing partitions}		\label{ncp}

The non-crossing partition lattice was first studied by
Kreweras~\cite{kre:pnc} who showed its M\"obius function to be a
Catalan number.  By using NBB sets we can combinatorially explain
this fact and relate these bases to the standard NBC bases for the
ordinary partition lattice.

If $\pi$ is a partition of $[n]:=\{1,2,\ldots,n\}$ into $k$ subsets,
or blocks, then we write $\pi=B_1/\ldots/B_k\ptn [n]$.  
When it will cause no confusion, we will not explicitly write out any
blocks that are singletons.  The set of
$\pi\ptn [n]$ form a lattice $\Pi_n$ under the refinement ordering.
We say that $\pi$ is {\sl non-crossing} if there do not exist $i,k\in B$ and
$j,l\in C$ for two distinct blocks $B,C$ of $\pi$ with $i<j<k<l$.  
Otherwise $\pi$ is {\sl crossing}.  The set of
non-crossing partitions of $[n]$ forms a meet-sublattice $NC_n$ of
$\Pi_n$ with the same rank function.  However unlike $\Pi_n$, $NC_n$
is not semimodular in general since if $\pi=13$ and $\si=24$ then
$\pi\mt\si=\zh$ and $\pi\jn\si=1234$ so
	$$\rho(\pi)+\rho(\si)=2<3=\rho(\pi\mt\si)+\rho(\pi\jn\si).$$
Here and throughout this paper, semimodularity refers to
upper-semimodularity. 

Another way to view non-crossing partitions will be
useful.  Let $G=(V,E)$ be a graph with vertex set $V=[n]$.  Then $G$ is 
{\sl non-crossing} if there do not exist edges $ik,jl\in E$ with
$i<j<k<l$.  Equivalently, $G$ is non-crossing if, when the vertices are
arranged in their natural order clockwise around a circle and the
edges are drawn as straight line segments, no two edges of $G$
cross geometrically.  Given a partition $\pi$ we can form a graph
$G_{\pi}$ by representing each block $B=\{i_1<i_2<\ldots<i_l\}$ by a
cycle with edges $i_1i_2,i_2i_3,\ldots,i_li_1$.  (If $|B|=1$ or 2 then
 $B$ is represented by an isolated vertex or edge, repectitively.)
Then it is easy to see that $\pi$ is non-crossing as a partition if and
only if $G_\pi$ is non-crossing as a graph.  In Figure~\ref{ptg} we
have displayed some partitions and their graphs.  

\bfi
\btab{rcc}
$\pi:$	&$1567/234/8$ (noncrossing)	&$134/256/78$ (crossing)\\
&
\bpi(70,70)(0,0)
\put(-35,35){\makebox(0,0){$G_\pi:$}}
\GcaL{-20pt}{6}
\GeaL{20pt}{5}
\GacL{-20pt}{7}
\GgcL{20pt}{4}
\GaeL{-20pt}{8}
\GgeL{20pt}{3}
\GcgL{-20pt}{1}
\GegL{20pt}{2}
\Gcaea \Gcaac \Geacg \Gaccg \Ggcge \Ggceg \Ggeeg
\epi
&
\bpi(110,70)(-20,0)
\GcaL{-20pt}{6}
\GeaL{20pt}{5}
\GacL{-20pt}{7}
\GgcL{20pt}{4}
\GaeL{-20pt}{8}
\GgeL{20pt}{3}
\GcgL{-20pt}{1}
\GegL{20pt}{2}
\Gcaea \Gcaeg \Geaeg \Gacae \Ggcge \Ggccg \Ggecg
\epi
\etab
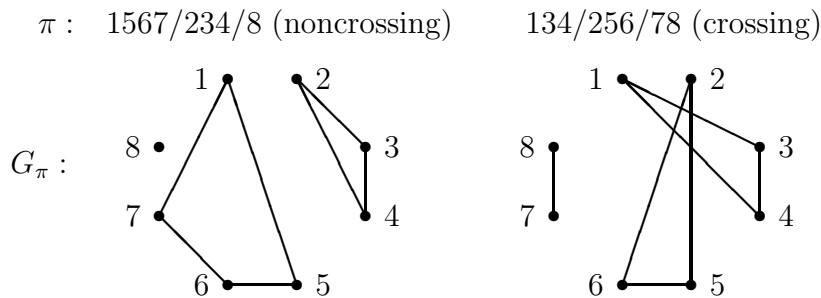
\capt{Partitions and their graphs}		\label{ptg}
\efi

The atoms of $NC_n$ are the partitions of the form $\de=ij$ where we
will always assume $i<j$.  Then $G_\de$ is a single edge, so we can
consider any $B\sbe A(NC_n)$ as a graph $G_B$ with an edge for each
$\de\in B$.  Define $ij\lt i'j'$ if and only if $j<j'$.
So the poset $(A(NC_n),\lte)$ will be ranked with the
elements at rank $j-2$ being all atoms of the form $ij$.  
An example for $n=4$ is given in Figure~\ref{poA}

\setlength{\unitlength}{3pt}
\bfi
\bpi(100,60)(-35,0)
\GaaLV{-35pt}{12}{1.5}
\GacLV{-35pt}{13}{1.5}
\GccLV{35pt}{23}{1.5}
\GaeLV{-35pt}{14}{1.5}
\GceLV{-30pt}{24}{1.5}
\GeeLV{35pt}{34}{1.5}
\Gaaac \Gaacc \Gacae \Gacce \Gacee \Gccae \Gccce \Gccee
\epi
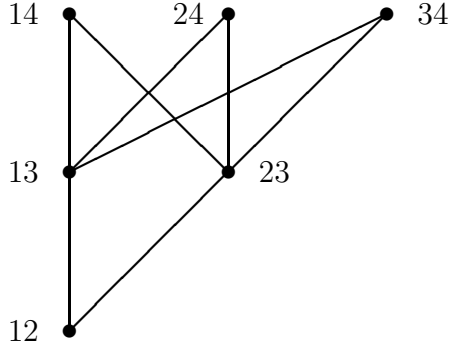
\capt{The partial order for $A(NC_4)$}		\label{poA}
\efi

In order to characterize the NBB sets, we first need a lemma.
\ble
Let $\de,\de'$ be a pair of atoms such that  either $\de$ and $\de'$
have the same rank 
in $(A(NC_n),\lte)$ or the graph of $D=\{\de,\de'\}$ is crossing.
Then $D$ is BB.
\ele
\pf\  Suppose first that the two atoms have the same rank.  Then
$\de=ij$ and $\de=i'j$ where without loss of generality $i<i'<j$.
Let $\al=ii'$.  Then $\de,\de'\gte\al$ and $\de\jn\de'=ii'j>\al$ so
$D$ is BB.

Suppose instead that the given atoms are crossed.  So we have $\de=ij$
and $\de=i'j'$ with $i<i'<j<j'$.  Letting $\al=ii'$ we get
the same inequalities as before, noting that $\de\jn\de'=ii'jj'$.\hfill\Qed

\bth							\label{NCNBB}
The NBB bases of $\oh$ in $NC_n$ are
all $B$ obtained by picking exactly one element from each rank of
$(A(NC_n),\lte)$ so that the corresponding graph $G_B$ is non-crossing.
\eth
\pf\  First suppose that $B$ is an NBB base of $\oh$.  Then by the
previous lemma we know that $B$ contains at most one element from each
rank and that $G_B$ is non-crossing.  If we do not pick an element
from some rank then $G_B$ is not connected.  But such a non-crossing
graph has a block of $\Jn B$ for each component of $G_B$,
contradicting $\Jn B=\oh$.

Conversely, suppose $B$ is picked according to the two given rules.
Then $G_B$ is connected and so
$\Jn B=\oh$.  If $B\spe D$ with $D$ a BB set we will derive a
contradiction.  Let $\de\in D$ be $\lte$-minimal and let $\al=ij$ be the
corresponding element guaranteed by the definition of bounded below.
Then~(\ref{BB1}) and our choice of $\de$ shows that for any
$\de'=i'j'\in D$ we have $i,j<j'$.  Since $B$ is non-crossing, so is
$D$ and thus $\Jn D$ is the same in $NC_n$ and $\Pi_n$.  It follows
from~(\ref{BB2}) that there is a path in $G_D$ of the form
$i=i_0,i_1,\ldots, i_l=j$ where each edge $\{i_ki_{k+1}\}$ is an atom of
$D$.  But our remarks about $\de'$ imply that $i_0<i_1$ and
$i_{l-1}>i_l$ so there must be an index $m$ such that
$i_{m-1}<i_m>i_{m+1}$.  Thus $D$ has two elements from the same
$\lte$-rank, the promised contradiction.\hfill\Qed

Note that the graphs $G_B$ in the previous theorem are certain
spanning trees on the vertex set $[n]$.  Furthermore, to get the NBB
bases for all elements 
of $NC_n$ it suffices to use the non-crossing restriction but picking at most
one element from each rank.  Finally, if one removes the non-crossing
restriction one gets exactly the standard NBC bases for the geometric
lattice $\Pi_n$. 

It is now easy to compute the M\"obius function of $NC_n$.  
It suffices to do this for $\oh$ since for any 
$\pi=B_1/\ldots/B_k\in NC_n$, the interval $[\zh,\pi]\iso \prod_i
NC_{|B_i|}$. 
Recall that the {\sl Catalan numbers} are defined by 
	$$C_n=\frac{1}{n+1}{2n\choose n}.$$
\bco[Kreweras]
We have
	$$\mu(NC_n)=(-1)^{n-1}C_{n-1}.$$
\eco
\pf\  All trees on $n$ vertices have $n-1$ edges.  Furthermore it is
easy to see that the number $T_n$ of non-crossing trees on $[n]$ of
the given form  and
$C_{n-1}$ satisfy the same initial conditions and recurrence relation,
	$$T_n=\case{1}{if $n=1$}{\sum_{0<i<n} T_i T_{n-i}}{if $n>1$.}$$
The result now follows from Theorems~\ref{NBB} and~\ref{NCNBB}.\hfill\Qed

There is another natural ordering of $A(NC_n)$ for which a result
similar to 
Theorem~\ref{NCNBB} holds, namely $ij\lte i'j'$ if and only if
$[i,j]\spe[i',j']$ as intervals of integers.  With this ordering, NBB
bases of $\oh$ again correspond to trees with vertex set $[n]$, but
this time they are the non-crossing trees in which each vertex is
either greater 
than all its neighbors or less than all its neighbors.  It follows
easily that, in any such tree, $\{1,n\}$ is an edge, and deletion of
this edge leaves two smaller such trees with vertex sets $[k]$ and
$[k+1,n]$ for some $k$.  This observation easily implies that the
numbers of such trees satisfy the recurrence for the Catalan numbers.
In fact, there is a simple bijection between these
trees and proper parenthesizations $P$ of the 
product of $n$ factors (one of the most familiar interpretations of the 
Catalan numbers).  To describe the bijection, identify the $n$ 
factors with the elements of $[n]$, in order.  For any particular 
parenthesization $P$, consider the sub-products defined by $P$.  For 
example, if $n=5$ and $P$ is $(((12)3)(45))$, then the sub-products
are $(12), ((12)3), (45)$, and $(((12)3)(45))$ itself.  Now to build
the corresponding $T$, take the vertex set $[n]$ and draw, for each 
sub-product of $P$, an edge from the first to the last element of the 
subproduct.  In the previous example, the edges would be $12, 13, 45$,
and $15$.  
It is a variant of this ordering which allows us to compute $\mu$ for
the non-crossing $B_n$ and $D_n$ lattices.

\section{Non-crossing $B_n$ and $D_n$}		\label{ncb}

In this section, we apply Theorem~\ref{NBB} to calculate the M\"obius
invariants $\mu(\oh)$ of the non-crossing $B_n$ and $D_n$ (and
intermediate) lattices.  These lattices were introduced by
Reiner~\cite{rei:npc} who computed $\mu$ using generating functions.
Non-crossing $B_n$ consists of those partitions $\pi$ of
$\{1,2,\dots,n,-1,-2,\dots,-n\}$ that satisfy three conditions: First,
$\pi$ is invariant under the involution $k\mapsto-k$.  Second, at most
one block of $\pi$ is fixed by this involution; if there is such a
block, it is called the {\sl zero-block\/} of $\pi$.  Third, the
partition is non-crossing (as in Section~\ref{ncp}), with respect to
the ordering $1<2<\dots<n<-1<-2<\dots<-n$.  (The first two of these
conditions determine a lattice isomorphic to that associated with the
hyperplane arrangement $B_n$.)  Non-crossing $D_n$ is the subposet
consisting of those $\pi$ for which the zero-block, if present, does
not consist of only a single pair $\{k,-k\}$.  An intermediate lattice
can be associated to every subset $S\subseteq[n]$, by allowing the
zero-block to be $\{k,-k\}$ only if $k\notin S$.  We use the
notations $NCB_n$, $NCD_n$, and $NCBD_n(S)$ for these lattices; thus
$NCB_n=NCBD_n(\emptyset)$ and $NCD_n=NCBD_n([n])$.

We begin by calculating $\mu(\oh)$ for $NCB_n$ and afterward indicate
the minor changes needed to handle the rest of these lattices.  The
atoms of $NCB_n$ are of three sorts.  First, there are the partitions
where one block consists of two positive numbers, say $\{i,j\}$,
another block is $\{-i,-j\}$, and all the remaining blocks are
singletons.  Second, there are the partitions where one block consists
of a positive and a negative number, say $\{i,-j\}$, another is
$\{-i,j\}$, and the rest are singletons.  Third, there are the
partitions whose only non-singleton block is of the form $\{i,-i\}$.
Following Zaslavsky~\cite{zas:grs}, 
we depict atoms as signed edges in
a graph whose vertex set is $[n]$.  An atom of the first sort is
depicted as a positive edge $ij$, one of the second sort is depicted
as a negative edge $ij$, and one of the third sort is depicted as a
(negative) half-edge at $i$.  (There is no such thing as a positive
half-edge.)  To avoid confusion, we emphasize that these signed graphs
are quite different from the graphs with vertex set
$\{1,2,\dots,n,-1,-2,\dots,-n\}$ used in deciding whether a partition
is non-crossing.  We sometimes identify an atom with the corresponding
edge; in particular, we may refer to an atom as positive or negative.

We partially order the atoms as follows.  Associate to each atom,
depicted as a signed edge $ij$, the interval $[i,j]\subseteq [n]$; in
the case of a half-edge at $i$, the interval consists of just $i$.
Then define $a\lt b$ to mean that either the interval associated to
atom $a$ properly includes that associated to atom $b$ or the two
intervals are equal and $a$ is negative while $b$ is positive.
(Notice that, apart from signs and half-edges, this matches the
ordering described at the end of Section~\ref{ncp}.)

Regarding atoms as signed edges, we regard sets of atoms as signed
graphs with vertex set $[n]$.  It is not difficult (though a bit
tedious) to verify that every NBB set $B$ has, as a graph, the
following properties. 
\ben
\item[(i)] No two signed edges cross.
\item[(ii)] No vertex $i$ has a neighbor $j<i$ and also a neighbor
$j'>i$.  
\item[(iii)] Any path joining two vertices $i$ and $j$ and consisting
entirely of vertices $k$ with $i\leq k\leq j$ must be
just a single signed edge $ij$.
\item[(iv)] There is no cycle of length $\geq3$.  (A 2-cycle,
consisting of a positive edge and a negative edge in the same place,
is permitted; we refer to such a pair of signed edges as a double
edge.) 
\item[(v)] No two negative atoms in $B$ have disjoint associated
intervals. 
\item[(vi)] If the interval associated to one atom in $B$ is properly
included in the interval associated to another, and if the former
atom is negative, then so is the latter.
\item[(vii)] If there is a negative atom in $B$, then there is exactly
one, say $a$, that is $\lt$-maximal (i.e., its interval is
inclusion-minimal); all other atoms in $B$ are negative if their
intervals properly include that of $a$ and positive otherwise.
\item[(viii)] $B$ has at most one half-edge, has at most one
double-edge, and cannot have both.
\een
Actually, only items (i), (ii), (v), and (vi) in this list directly
use the NBB assumption.  The other four items follow from these purely
graph-theoretically.

Conversely, any signed graph satisfying (i) through (viii) is NBB when
viewed as a set of atoms of $NCB_n$.  We leave the verification to the
reader, with the hint that items (i), (v), and (vi) ensure that the
join of any atoms from this set is the same whether computed in
$NCB_n$ or in $B_n$, because none of these atoms cross when regarded
as partitions of $\{1,2,\dots,n,-1,-2,\dots,-n\}$.

{}From this characterization of the NBB sets in $NCB_n$, we easily
obtain a characterization of the NBB bases of $\oh$.  These bases
(regarded as signed graphs) are obtainable as follows.  First, take a
non-crossing tree $T$ with vertex set $[n]$ in which each vertex is
either greater than all its neighbors or less than all its neighbors.
(This part is just as at the end of Section~\ref{ncp}.)  Then pick
either an edge or a vertex of $T$.  If you picked an edge $e$, then
make it a double edge, give negative signs to all the edges of $T$
whose intervals properly include that of $e$, and give all remaining
edges of $T$ positive signs.  If you picked a vertex $v$, then attach
a (negative) half-edge at $v$, give negative signs to all edges of $T$
whose interval contains $v$, and give all remaining edges of $T$
positive signs.

Finally, to apply Theorem~\ref{NBB}, we count the NBB bases for $\oh$,
i.e., we count the signed trees of the sort just described.  We
already saw in Section~\ref{ncp} that there are 
$$
C_{n-1}=\frac{1}{n}{2n-2\choose n-1}=\frac{1}{2n-1}{2n-1\choose n}
$$
ways to choose $T$.  Then there are $2n-1$ ways to choose an edge or
vertex, since there are $n$ vertices and $n-1$ edges.  After this
choice, the rest of the construction of the NBB signed graph is
completely determined.  So the number of NBB bases for $\oh$ is 
$$
C_{n-1}\cdot(2n-1)={2n-1\choose n}.
$$
Every NBB base of $\oh$ has exactly $n$ elements, namely the
$n-1$ edges of $T$ (with signs) plus either an extra edge if you chose
an edge and doubled it or an extra half-edge if you chose a vertex.
Therefore, $\mu(NCB_n)=(-1)^n{2n-1\choose n}$.

The calculation for $NCBD_n(S)$ is almost exactly the same.  The only
difference is that half-edges can occur only at vertices not in $S$.
Thus, in the description of NBB bases for $\oh$, we need only replace 
``pick an edge or a vertex of $T$'' with ``pick an edge of $T$ or a
vertex not in $S$.''  Thus, the number of options at this step is no
longer $2n-1$ but only $2n-1-|S|$.  Therefore, we obtain, in agreement
with Reiner's calculation (\cite{rei:npc}), the M\"obius invariant
for $NCBD_n(S)$:
$$
\mu(NCBD_n(S))=(-1)^nC_{n-1}\cdot(2n-1-|S|).
$$

\section{Shuffle Posets}			\label{sp}

The poset of shuffles was introduced by Greene~\cite{gre:ps}.  We
need to recall some of his definitions and results before applying
Theorem~\ref{NBB}.  Let $\cA$ be a set, called the {\sl alphabet of
letters}.  A 
{\sl word over $\cA$} is a sequence $\bu=u_1 u_2\ldots u_n$
of elements of $\cA$.  All of our words will
consist of distinct letters and we will sometimes also
use $\bu$ to stand for the set of letters in the word, depending upon
the context.  
A {\sl subword of $\bu$} is $\bv=u_{i_1}\ldots u_{i_l}$ where
$i_1<\ldots<i_l$. If $\bu,\bv$ are any two words then the {\sl
restriction of $\bu$ to $\bv$} is the subword $\bu_{\bv}$ of $\bu$
whose letters are exactly those of $\bu\cap\bv$.  
A {\sl shuffle of $\bu$ and $\bv$} is any word $\bs$ such that
$\bs=\bu\uplus\bv$ as sets (disjoint union) and $\bs_{\bu}=\bu$,
$\bs_{\bv}=\bv$.

Given nonnegative integers $m,n$ Greene defined the {\sl poset of shuffles}
$\cW_{m,n}$ as follows.  Fix disjoint words $\bx=x_1\ldots x_m$ and 
$\by=y_1\ldots y_n$.  The elements of $\cW_{m,n}$ are all shuffles
$\bw$ of a subword of $\bx$ with a subword of $\by$.  The partial
order is $\bv\le\bw$ if $\bv_{\bx}\spe\bw_{\bx}$ and
$\bv_{\by}\sbe\bw_{\by}$.  It is easy to see that $\cW$ has minimal
element $\zh=\bx$, maximal element $\oh=\by$ and is ranked with rank
function
	\beq						\label{rhoW}
	\rho(\bw)=(m-|\bw_{\bx}|)+|\bw_{\by}|.
	\eeq
For example, $\cW_{2,1}$ is shown in Figure~\ref{W}
where $\bx=de$ and $\by=D$.

\setlength{\unitlength}{3pt}
\bfi
\bpi(50,70)(0,-10)
\GdaV{1.5}
\put(30,-5){\makebox(0,0){$de$}}
\GacLV{-15pt}{$d$}{1.5}
\GbcLV{-15pt}{$e$}{1.5}
\GccLV{30pt}{$Dde$}{1.5}
\GecLV{30pt}{$dDe$}{1.5}
\GgcLV{30pt}{$deD$}{1.5}
\GaeLV{-15pt}{$\emp$}{1.5}
\GceLV{-25pt}{$Dd$}{1.5}
\GdeLV{25pt}{$De$}{1.5}
\GfeLV{-25pt}{$dD$}{1.5}
\GgeLV{25pt}{$eD$}{1.5}
\GdgV{1.5}
\put(30,65){\makebox(0,0){$D$}}
\Gdaac \Gdabc \Gdacc \Gdaec \Gdagc
\Gacae \Gacce \Gacfe \Gbcae \Gbcde \Gbcge
\Gccce \Gccde \Gecde \Gecfe \Ggcfe \Ggcge
\Gaedg \Gcedg \Gdedg \Gfedg \Ggedg
\epi
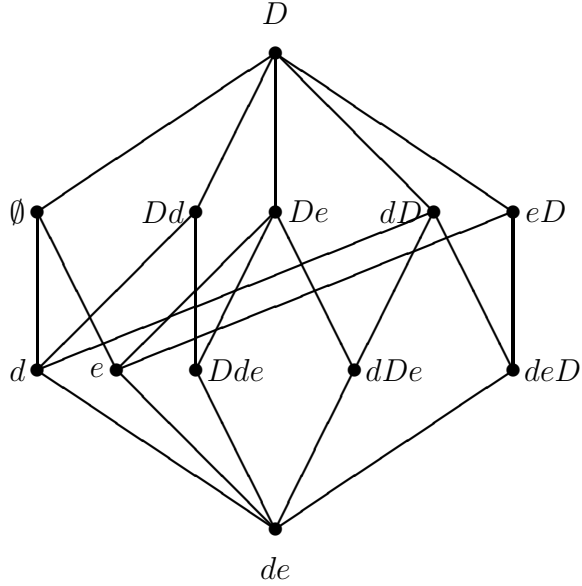
\capt{The lattice $\cW_{2,1}$} \label{W}
\efi

In order to apply the NBB Theorem we will need to describe the join
operation in $\cW_{m,n}$.  Greene does this using crossed elements,
not to be confused with the crossing partitions discussed in 
Section~\ref{ncp}.
Given $\bu,\bv\in\cW_{m,n}$ then $x\in\bx$ is {\sl crossed
in $\bu,\bv$} if there exist $y_i,y_j\in\by$ such that $i\leq j$ and
$x$ appears before $y_i$ in one of the two words but after $y_j$ in the other.
For example if $\bx=def$ and $\by=DEF$ then for the pair
$\bu=dDEe,\bv=Fdef$ the only crossed letter is $d$.  The join of
$\bu,\bv$ is then the unique word $\bw$ greater than both $\bu,\bv$ such that
	\beq					\label{Wjn}
	\barr{rcl}
	\bw_{\bx}&=&\{x\in\bu_{\bx}\cap\bv_{\bx}\ :\ 
		x \mbox{ is not crossed}\}\\
	\bw_{\by}&=&\bu_{\by}\cup\bv_{\by}.
	\earr
	\eeq
In the previous example, $\bu\jn\bv=DEFe$.  This example also shows
that $\cW_{m,n}$ is not geometric since the semimodularity law is violated:
	$$\rho(\bu)+\rho(\bv)=3+1<5=\rho(\bu\jn\bv)\le
	\rho(\bu\jn\bv)+\rho(\bu\mt\bv).$$

The set atoms $A_{m,n}=A(\cW_{m,n})$ consists of two types.  An {\sl
a-atom}, 
respectively {\sl b-atom}, is one obtained from $\bx$ by deleting a
letter of $\bx$, respectively inserting a letter of $\by$.  Let $A_a$
denote the set of a-atoms and similarly for $A_b$.
Define $\lte$ on $A_{m,n}$ to be the poset whose relations are all
those of the form $\ba\lt\bb$ with $\ba\in A_a$ and  $\bb\in A_b$ 
\ble							\label{WBB}
Suppose $\bb,\bb'\in A_b$.  If $\bb,\bb'$ have crossed elements then
$D=\{\bb,\bb'\}$ is a BB set of $(A_{m,n},\lte)$.
\ele
\pf\  Our hypothesis and~(\ref{Wjn}) show that $(\bb\jn\bb')_{\bx}\sbs\bx$
(proper containment of sets).  So there is an a-atom $\ba$ with
$\ba\spe (\bb\jn\bb')_{\bx}$  which forces $\ba<\bb\jn\bb'$ in $W_{m,n}$.
Since by definition $\ba\lt\bb,\bb'$ the element $\ba$ satisfies the
definition of a BB set.\hfill\Qed

The next result will characterize the NBB bases and show that the 
converse of the previous lemma also holds.
\bth							\label{WNBB}
Let $\bs$ be a shuffle of $\bx,\by$ and consider
	$$B_{\bs}=A_a\cup\{\bb\in A_b\ :\ \bb\leq\bs\}.$$
Then the NBB bases of $\by\in\cW_{m,n}$ under the given partial order
are exactly the $B_{\bs}.$
\eth
\pf\  Suppose first that $B$ is an NBB base of $\by$.  Then for each
element $y\in\by$ we must have a corresponding b-atom $\bb_y$ in order
to get $\Jn B=\by$.  In fact there must be exactly one such atom for
each $y\in\by$ and these atoms must all lie below a shuffle $\bs$, for
otherwise $B$ would contain a BB pair as in Lemma~\ref{WBB}.  It 
follows that $\Jn_y \bb_y=\bs$.  So in order to get $\Jn B=\by$ we
must have $A_a\sbe B$.  Thus $B$ is of the form $B_{\bs}$ as desired.

Conversely, consider any $B_{\bs}$.  It is easy to see that 
$\Jn B_{\bs}=\by$. We will show that $B_{\bs}$ is NBB by contradiction.
Suppose that $D\sbe B_{\bs}$ is a BB set.  Then~(\ref{BB1}) forces $D$
to contain only b-atoms which in turn implies $\Jn D\le\bs$.  Now pick any
$d\in D$ and let $a\in A_{m,n}$ be the atom guaranteed by the
definition of BB.  Then $a$ is an a-atom and $a<\Jn D\le\bs$.  But
this contradicts the fact that there are no a-atoms below any shuffle of
$\bx$ and $\by$.  Thus $B_{\bs}$ is an NBB base of $\by$ and our
characterization is complete.\hfill\Qed

To determine the M\"obius function of $\cW_{m,n}$ it suffices to
compute $\mu(\oh)$ since for any $\bw\in\cW_{m,n}$ the interval
$[\zh,\bw]$ is isomorphic to a product of $\cW_{p,q}$'s for 
certain $p\le m$ and $q\le n$.
\bco[Greene]
We have
	$$\mu(\cW_{m,n})=(-1)^{m+n}{m+n\choose m}.$$
\eco
\pf\  The bases $B_{\bs}$ all have size $m+n$ and the number of
possible shuffles $\bs$ is ${m+n\choose m}$.  The result now follows
from Theorems~\ref{NBB} and~\ref{WNBB}.\hfill\Qed

\section{The Dominance Order}	\label{dom}

Bogart~\cite{bog:mfd} and Brylawski~\cite{bry:lip} first computed the
two-variable M\"obius function of the lattice of integer partitions under
dominance.  Subsequently Greene~\cite{gre:clm} gave two alternative
ways to compute this function.  Our NBB set characterization leads to
a formula in Corollary~\ref{Pmu} which is essentially equivalent to, but
simpler than, Greene's second description of 
$\mu$~\cite[Theorem 4.1]{gre:clm}.

We begin by reviewing the relevant
definitions.  A {\sl partition\/} $\la$ of $n$
is a weakly decreasing sequence of positive integers
$\la_1\geq\la_2\geq\dots\geq\la_r$ whose sum is $n$.  For any such
partition $\la$ and any non-negative integer $k$, we write $|\la|_k$
for $\sum_{i=1}^k\la_i$, where $\la_i$ is interpreted as 0 for $i>r$,
so $|\la|_k=n$ for all $k\geq r$.  A partition $\la$ {\sl
dominates\/} another partition $\nu$ (of the same $n$), written
$\la\geq\nu$, if $|\la|_k\geq|\nu|_k$ for all $k$.  In terms of
Ferrers diagrams (in the English orientation) $\la\geq\nu$ means
that the diagram for $\la$ can be obtained from that of $\nu$ by
moving some squares up to earlier rows.

It is well known that this ordering makes the set of partitions of $n$
into a lattice $\cP_n$; we review the construction of joins in $\cP_n$
because it will be needed for our NBB calculations.  A {\sl
composition\/} of $n$ is like a partition except that the parts
$\la_i$ need not be in weakly decreasing order.  The dominance order of
compositions is defined exactly as for partitions, and the result is a
lattice in which joins are easily computed since
$|\la\jn\nu|_k=\max\{|\la|_k,|\nu|_k\}$.  The join of two
partitions, computed in this lattice of compositions, need not be a
partition.  However, for every composition $\la$ there is a unique
smallest partition above $\la$, which we call the {\sl partition
reflection\/} of $\la$.  Joins of partitions in $\cP_n$ can be
computed by first forming the join in the dominance lattice of
compositions and then forming the partition reflection of the result.
We should point out that the partition reflection
of a composition need not be simply the result of rearranging the
parts into weakly decreasing order.  For example, the partition
reflection of $(1,3)$ is not $(3,1)$ but $(2,2)$.

The bottom element $\zh$ of $\cP_n$ is the partition $(1,1,\dots,1)$
and the only atom is $(2,1,1,\dots,1)$.  By Theorem~\ref{cct},
these are the only elements of $\cP_n$ where the M\"obius function has
a non-zero value. But we can also consider the M\"obius function of
elements in any upper interval $[\be,\oh]$; 
indeed, the value of this M\"obius
function at some $\la\geq\be$ is what is usually called the
(two-variable) M\"obius function $\mu(\be,\la)$.

To calculate this M\"obius function, we fix $\be$, we describe the
atoms of $[\be,\oh]$, and we calculate the joins of sets of atoms with
particular attention to determining when one atom is below a join of
others.  Then, we describe an appropriate partial ordering $\lte$ of
the atoms, characterize its NBB sets, and use Theorem~\ref{NBB} to
evaluate the M\"obius function.  

For the rest of this section, let $\be=(\be_1,\be_2,\dots,\be_r)$ be a
fixed but arbitrary partition of $n$.  By a {\sl wall\/} we mean a
maximal sequence of at least two numbers $k$ for which the
corresponding $\be_k$ are equal.  (This corresponds to a {\sl flat} in
Greene's terminology.)
Thus, an interval
$[i,j]\subseteq[1,r]$ is a wall if and only if $i<j$ and 
$\be_{i-1}>\be_i=\be_j>\be_{j+1}$.  Here and below, we use the
convention that $\be_0=\infty$ and $\be_{r+1}=0$, so that 1
(respectively, $r$) can be part of a wall if the first (respectively, last) two
components of $\be$ are equal.  If $[i,j]$ is a wall then we call $i$ its
top and $j$ its bottom, the terminology being suggested by the Ferrers
diagram.

The atoms of $[\be,\oh]$ are of two sorts: 
\ben 
\item[(i)] If $1\leq i<r$ and
neither $i$ nor $i+1$ is in a wall (i.e.,
$\be_{i-1}>\be_i>\be_{i+1}>\be_{i+2}$), then $\be$ is covered by the
partition $\al$ that agrees with $\be$ except that $\al_i=\be_i+1$ and
$\al_{i+1}=\be_{i+1}-1$.  We denote this atom $\al$ of $[\be,\oh]$ by
$i+1\ra i$, since its Ferrers diagram is obtained from that of $\be$
by moving a square from row $i+1$ up to row $i$.  
\item[(ii)] If $[i,j]$ is a
wall, then $\be$ is covered by the partition that agrees with $\be$
except that $\al_i=\be_i+1$ and $\al_j=\be_j-1$.  We denote this $\al$
by $j\ra i$.  
\een

For both types of $\al$, we refer to the set of $k$ where
$|\al|_k\neq|\be|_k$ as the {\sl critical interval\/} of $\al$, $I_\al$.
(This is {\sl not\/} the same as Brylawski's use of ``critical.'')
It consists of only $i$ for $\al$ of type (i) and it is $[i,j-1]$ for
$\al$ of type (ii).  Note that $|\al|_k=|\be|_k+1$ if and only if $k\in I_\al$
Note also that different atoms of $[\be,\oh]$
have disjoint critical intervals and
	\beq						\label{Ial}
	\biguplus_\al I_\al = [1,r-1]\setm
	\{k:\ \mbox{$[i,k]$ or $[k+1,i]$ is a wall for some $i$}\}.
	\eeq
So an element of $[1,r-1]$ lies in a critical
interval unless it is the bottom of a wall or the immediate
predecessor of the top of a wall.

Let us consider the join of some subset $\cB\sbe A([\be,\oh])$, first
in the sense of compositions and then in the sense of partitions.  
If $\ga$ is the composition join, then from the previous paragraph
	\beq						\label{ga}
	|\ga|_k - |\be|_k =\case{0}{if $k\not\in X$}{1}{if $k\in X$}
	\eeq
where $X=\uplus_{\al\in\cB} I_\al$.  If this $\ga$ is a partition,
then it is also the join in the partition sense; 
otherwise, we must take its partition reflection.  So we consider next
how $\ga$ could fail to be a partition and how it is changed by
reflection.

By~(\ref{ga}) the components of $\ga$ differ from those of $\be$ by $\pm1$
or 0.  It follows that any failure of $\ga$ to be a partition, i.e.,
any occurrence of $\ga_k<\ga_{k+1}$, must arise in one of two ways:
$\be_k=\be_{k+1}$ or $\be_{k+1}+1$.
The first possibility can be excluded, because it requires $k$ and
$k+1$ to be part of a wall, say $[i,j]$.  
If the critical  interval $[i,j-1]\sbe X$ then 
by~(\ref{Ial}) and~(\ref{ga}) we have $\ga_i=\be_i+1$,
$\ga_j=\be_j-1$, and $\ga$ agrees with $\be$ on the rest of $[i,j]$;
if, on the other hand, $[i,j-1]\cap X=\emp$, then $\ga$
agrees with $\be$ on all of $[i,j]$.  
In either case, we cannot have $\ga_k<\ga_{k+1}$.
So any occurrence of $\ga_k<\ga_{k+1}$ has $\be_k=\be_{k+1}+1$,
$\ga_k=\be_k-1=\be_{k+1}$, and $\ga_{k+1}=\be_{k+1}+1=\be_k$.  Thus
by~(\ref{ga}) $X$ must contain $k-1$ and $k+1$ but not $k$.

Consider the set $Y$ obtained by adjoining to $X$ all those numbers
$k\notin X$ for which $k-1,k+1\in X$ and $\be_k=\be_{k+1}+1$.  Let
$\nu$ be the composition of $n$ such that
	\beq						\label{nu}
	|\nu|_k - |\be|_k =\case{0}{if $k\not\in Y$}{1}{if $k\in Y$.}
	\eeq
The preceding
observations show that $\nu$ is a partition, because, by enlarging $X$
to $Y$, we have corrected all the failures of $\ga$ to be a partition.
Furthermore, $\nu$ clearly dominates $\ga$.  We claim that $\nu$ is
the partition reflection of $\ga$, i.e., that every partition $\eta$
dominating $\ga$ also dominates $\nu$.  To see this, it suffices to
show that $|\eta|_k>|\ga|_k$ for $k\in Y-X$; but if this failed for
some $k$, then $\eta$ would fail to be a partition because
$\eta_k<\eta_{k+1}$ just as for $\ga$.  Thus, $\nu$ is the join of
$\cB$ in $\cP_n$.

We must also determine which 
atoms of $[\be,\oh]$ other than members of $\cB$ are below $\nu$.  
Distinct atoms have disjoint critical intervals,  so the atoms
we are looking for are those whose critical intervals are singletons
$k$ where $k\in Y-X$.  Note that there is an
atom with critical interval $k$ if and only if either $k$ and
$k+1$ constitute a wall or neither of them belongs to a wall.  The
former alternative is irrelevant in the present context, since $k\in
Y-X$ implies $\be_k=\be_{k+1}+1$.  So we need only consider the second
alternative, where $\be_{k-1}>\be_k>\be_{k+1}>\be_{k+2}$.  

We are ready to begin the description of the NBB bases for an
arbitrary $\la\geq\be$.  Let $\cA$ be the set of atoms of $[\be,\la]$.
If $\Jn\cA\neq\la$ then $\mu(\be,\la)=0$ by Theorem~\ref{cct}, so from
now on we assume 
$\Jn\cA=\la$.  Call an atom $b\in\cA$ of the form $k+1\ra k$ 
{\sl special\/}  
if $\be_k=\be_{k+1}+1$ and $\cA$ also contains atoms $a,c$ equal to 
$k+2\ra k+1,k\ra k-1$, respectively.  The preceding discussion shows
that $b<a\jn c$; 
conversely, if an atom in $\cA$ is under the join of some other
atoms in $\cA$, then it must be special.  We let $\cS\sbs\cA$ be the
set of special atoms.

Let us list the atoms in $\cA$ according to the ordering of their
critical intervals.  
By a {\sl special run\/} in this list, we mean a maximal sequence of
consecutive $\si\in\cS$.  Define a partial ordering $\lte$
on $\cA$ by imposing on each sequence
	\beq						\label{seq}
	\tau=\si_0,\si_1,\ldots,\si_{q+1}=\tau'
	\eeq
where $\si_1,\ldots,\si_q$ is  a special run (so $\tau,\tau'\not\in\cS$)
the relations
	$$\si_{3i+1}\lt\si_{3i},\si_{3i+2}\quad\mbox{for}\quad
				1\le3i+1\le q.$$
An example is given in Figure~\ref{run}.  
Note that this ordering makes $\{\si_{3i},\si_{3i+2}\}$ a BB set for
each $i$.
We let $\cS_j$, $0\le j\le2$,
be the set of elements in all special runs of the form $\si_{3i+j}$
for some $i$, so $\cS=\cS_0\uplus\cS_1\uplus\cS_2$.

\setlength{\unitlength}{3pt}
\bfi
\bpi(60,30)(0,-10)
\GabV{1.5} \GbaV{1.5} \GcbV{1.5} \GdbV{1.5} \GeaV{1.5} \GfbV{1.5} \GgbV{1.5}
\Gbaab \Gbacb \Geadb \Geafb
\put(0,15){\makebox(0,0){$\tau$}}
\put(10,-5){\makebox(0,0){$\si_1$}}
\put(20,15){\makebox(0,0){$\si_2$}}
\put(30,15){\makebox(0,0){$\si_3$}}
\put(40,-5){\makebox(0,0){$\si_4$}}
\put(50,15){\makebox(0,0){$\si_5$}}
\put(60,15){\makebox(0,0){$\tau'$}}
\epi
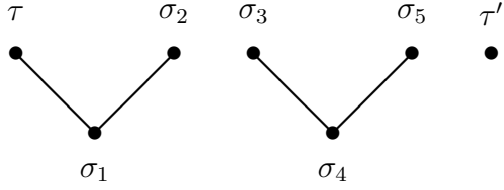
\capt{The partial order on an extended special run }	\label{run}
\efi

\bth							\label{PNBB}
Let $[\be,\la]$ be an interval in $\cP_n$ with atom set
$\cA$ such that $\Jn\cA=\la$.  
Then $\la$ has an NBB base in $\cA$ if and only if there is no special
run of length 1 modulo 3.  If a base exists, then it is unique and
equals 
	\beq						\label{B}
	B=\cA\setm\cS_2.
	\eeq
\eth
\pf\ If $B$ is an NBB base of $\la$ then it must include $\cA\setm\cS$
since these elements are not $\le$ the join of any others.  Now
consider any sequence of the form~(\ref{seq}).  Since $\tau\not\in\cS$
we have $\tau\in B$ and this forces $\si_2\not\in B$ since
$\{\tau,\si_2\}$ is BB.  Now we must put $\si_1,\si_3\in B$ since
neither is below $\Jn (\cA\setm\si_2)$.  Now repeat the argument
with $\si_3$ in place of $\tau$ to inductively see that the only
possible candidate for an NBB base is~(\ref{B}). But if there is a
special run~(\ref{seq}) of length $3k+1$ for some $k$ then $\si_{3k}$ and
$\tau'=\si_{3k+2}$ form a BB set in $B$, a contradiction.  If no
such run exists then we never have such a pair in $B$ and so it is NBB
as desired.\hfill\Qed

\bco							\label{Pmu}
Let $[\be,\la]$ be an interval in $\cP_n$ with atom set
$\cA$.  If $\Jn\cA\neq\la$ then $\mu(\be,\la)=0$.  If $\Jn\cA=\la$
then let $r_i$, $0\le i\le2$ be the number of special runs of length
$i$ modulo 3.  Then
	$$\mu(\be,\la)=
	\case{0}{if $r_1\ge1$}{(-1)^{|\cA\setm\cS|+r_2}}{if $r_1=0$.}
	$$
\eco
\pf\  We have already noted that $\mu(\be,\la)=0$ if $\Jn\cA\neq\la$.
For the other two cases, suppose first that $r_1\ge1$.
Then $\la$ has no NBB base in $[\be,\la]$ by Theorem~\ref{PNBB} so
$\mu(\be,\la)=0$ by Theorem~\ref{NBB}.  If $r_1=0$ then by the same
two results
	$$\mu(\be,\la)=(-1)^{|B|}=(-1)^{|\cA\setm\cS|+|\cS_0|+|\cS_1|}.$$
If a special run has length congruent to 0 (respectively, 2) modulo 3
then its contribution to $|\cS_0|+|\cS_1|$ is 0 (respectively, 1)
modulo 2. The last case now follows from the previous
displayed equation.\hfill\Qed

\section{Supersolvable Lattices} \label{sl}

Stanley \cite{sta:ssl} defined a {\sl supersolvable lattice\/} to be
a pair $(L,\Delta)$ where $L$ is a lattice,
$\Delta:\zh=x_0<x_1<\dots<x_{n-1}<x_n=\oh$ is a maximal chain of
$L$, and $\Delta$ together with any other chain of
$L$ generates 
a distributive lattice.  One often refers to $L$ as supersolvable,
$\Delta$ being tacitly understood.
It is easy to see that a supersolvable lattice has a rank function
$\rho$.
He showed that, if such an $L$ is
also semimodular, then its characteristic polynomial
\beq						\label{chi}
\chi(L,t)=\sum_{x\in L} \mu(\zh,x)t^{n-\rho(x)}
\eeq
factors as $(t-a_1)(t-a_2)\cdots(t-a_n)$, where $a_i$ is the number of
atoms of $L$ that are below $x_i$ but not below $x_{i-1}$.  
Our purpose in this section is to use Theorem~\ref{NBB} to prove
Stanley's factorization of the characteristic polynomial for 
a wider class of lattices.  For this purpose we will replace both
supersolvability and semimodularity by weaker hypotheses.

To state the first of our two hypotheses, we define an element $x$ of
a lattice $L$ to be {\sl left-modular\/} if, for all $y\leq z$,
$$
y\jn(x\mt z)=(y\jn x)\mt z.
$$
It is
standard to call $(x,z)$ a {\sl modular pair\/} if the preceding equation is
satisfied by every $y$ that is $\leq z$.  So $x$ is left-modular if
and only if every pair with $x$ on the left is modular.  Note that 
left-modularity, unlike the analogously defined right-modularity, is a
self-dual concept of lattice theory.  Our first hypothesis is that
\begin{center}
$L$ has a maximal chain $\Delta$, all of whose elements are
left-modular.  
\end{center}
In this case we will call $L$ itself {\sl left-modular}.  

Fix a maximal chain $\Delta:\zh=x_0<x_1<\dots<x_{n-1}<x_n=\oh$ 
whose elements may or may not be left-modular.  We partition the set
$A$ of atoms of $L$ into pieces 
$A_i=\{a\in A\mid a\leq x_i{\rm\ but\ }a\not\leq x_{i-1}\}$, which we
call the {\sl levels\/} of $A$, and we partially order $A$ by
setting $a\lt b$ if and only if $a$ is in an earlier level than $b$,
where ``earlier'' means having a smaller subscript.  
We say that this order $\lte$ is {\it induced\/} by the chain $\De$.
Note that an atom $a$
cannot be $\leq\Jn S$ for any set $S$ of atoms from strictly lower
levels, for there is an $x_i$ that is $\geq$ all elements of $S$ but
not $\geq a$.  We can now state our second
hypothesis.
\begin{center}
If $\lte$ is induced by $\De$ and $a\lt b_1\lt b_2\lt\dots\lt b_k$
then $a\not\leq\Jn_{i=1}^kb_i$.  
\end{center}
A lattice $L$ having a chain $\De$  with this property will be said to
satisfy the {\sl level condition}.  If $L$ is a finite lattice with a
maximal chain satisfying both the
left-modular 
and level conditions  then it is called an {\sl LL} lattice.

\bpr \label{usmb}
\ben
\item If $L$ is supersolvable then it is left-modular but not conversely.
\item If $L$ is semimodular then it satisfies the level condition
(for any maximal chain) but not conversely.
\een
\epr
\pf\ 
Stanley has already pointed out that first statement
holds~\cite[Proposition 2.2 and ff.]{sta:ssl}.
For the second,
recall that semimodularity implies that 
if $x\in L$ and $a$ is an atom not below $x$ then
$x\jn a$ covers $x$.  
Now prove the implication by induction on $k$, the number of $b$ atoms
in the level condition.  The 
cases $k=0$ and $k=1$ are obvious, so suppose the result holds for
$k-1$ but fails for $k$.  So we have $a\lt b_1\lt b_2\lt\dots\lt b_k$
and $a\leq\Jn_{i=1}^kb_i$.  Let $c=\Jn_{i=1}^{k-1}b_i$ and notice that
$a\not\leq c$ by induction hypothesis.  Thus $c<c\jn a\leq c\jn b_k$.
But $b_k$ is an atom, so $c\jn b_k$ covers $c$.  Therefore, $c\jn
a=c\jn b_k\geq b_k$.  But this is absurd, as $c\jn a$ is a join of
atoms from levels strictly earlier than that of $b_k$. 

The implication is not reversible as the shuffle
posets $\cW_{m,1}$ serve as counterexamples.  
It is not hard to check that the level condition holds
when $\by=y_1$ with $\Delta$ being the chain where the $x_i$'s are
first removed one at a time and then $y_1$ is added.  On the
other hand, as long as  
$m\geq1$, there are pairs of atoms not covered by their
join (atoms that add $y_1$ in different places), so $\cW_{m,1}$
is not a semimodular lattice.\hfill\Qed

We intend to generalize to LL lattices Stanley's factorization of the
characteristic polynomial.  Since LL lattices need not be ranked, we
need a suitable substitute for the rank function used in defining the
characteristic polynomial.  Let $L$ be a left-modular lattice, with
left-modular maximal chain $\De$ inducing a partition into levels
$A_i$ as above.
Define the {\sl generalized rank\/} of $x\in L$ to be
$$
\rho(x)=\mbox{ number of $A_i$'s containing atoms $\le x$.}
$$
Later in this section, we shall relate this $\rho$ to lengths
of chains, but for now we use it in (\ref{chi}) (with $n$ still
denoting the length of the chain $\De$) to define the characteristic
polynomial of any left-modular lattice.  We shall obtain a
factorization of this polynomial for any LL lattice.

We first characterize the NBB sets using the following
lemma. 
\ble     \label{2bb}
If $a$ and $b$ are distinct atoms from the same level $A_i$ in a
left-modular lattice, then
$a\jn b$ is above some atom $c$ from an earlier level.  
\ele
\pf\ 
Since $b$ is below $x_i$ but not below $x_{i-1}$, we have
$x_{i-1}<x_{i-1}\jn b\leq x_i$.  By maximality of the chain $\Delta$,
it follows that $x_{i-1}\jn b=x_i\geq a\jn b$.  Applying
left-modularity of $x_{i-1}$ with $y=b$ and $z=a\jn b$, we find that
$$
b\jn(x_{i-1}\mt(a\jn b))=(b\jn x_{i-1})\mt(a\jn b)=a\jn b.
$$
But $b\not\geq a$ so the right side of this equation is not $b$.
Since the left side is also not $b$, we have that
$x_{i-1}\mt(a\jn b)$ is not $\zh$ and so is above some atom $c$
satisfying the lemma.\hfill\Qed

\bth  						\label{ptv}
In an LL lattice, the NBB sets are exactly those subsets of $A$ that
have at most one member in each level $A_i$.
\eth
\pf\ 
The level condition immediately implies that, if $B\subseteq A$ has no two
members from the same level then no subset of it is BB.  Conversely,
if $B$ has two members at the same level then, by the lemma just
proved, those two constitute a BB set.\hfill\Qed

The following lemma is useful not only for our primary goal, factoring
the characteristic polynomial of an LL lattice, but also for relating
our generalized rank functions to lengths of chains.  
\ble  \label{ind}
Let $B$ be an NBB set in an LL lattice.
Then every atom
$a\leq\Jn B$ is at the same level as some element of $B$.  
In particular, any NBB base for $x$ has exactly $\rho(x)$
elements. 
\ele
\pf\ 
It suffices to prove the first statement, since the second follows
from it and Theorem~\ref{ptv}.  
So suppose $B$ and $a$ were a
counterexample to the first statement.  
Let $A_j$ be the level containing $a$, and let $y$ be the join of all
the elements of $B$ of higher level than $A_j$.  Since $B$ has no
element at level $A_j$, we have $a\leq\Jn B\leq x_{j-1}\jn y$.
Setting $z=a\jn y$, we obtain from the left-modularity of $x_{j-1}$
that $(y\jn x_{j-1})\mt z=y\jn(x_{j-1}\mt z)$.  

On the left side of this last equation, both sides of the meet are
$\geq a$, and therefore so is the whole left side.  On the right side, since
$z$ is a join of atoms from levels $A_j$ and higher,
the level condition tells us that $z$
is above no atom of lower level than $A_j$, so no atom is below
$x_{j-1}\mt z$.  Therefore $x_{j-1}\mt z=\zh$ and the right side
reduces to $y$.  

Combining these results, we have $a\leq y$.  But, in view of the
definition of $y$, this contradicts the level condition.\hfill\Qed

We can now prove our generalization
of Stanley's theorem on semimodular supersolvable lattices.
\bth						\label{LL}
If $L$ is an LL lattice then
characteristic polynomial of $L$ factors as
$$
\chi(L,t)=\prod_{i=1}^n(t-|A_i|)
$$
\eth
\pf\ 
Combining Theorems~\ref{NBB} and~\ref{ptv} along with the new
definition of $\rho$ we have
\bea
\chi(L,t)
&=&\sum_{x\in L}\ \sum_{{B{\rm\ NBB}}\atop{\Jn B=x}} (-1)^{|B|}t^{n-\rho(x)}\\
&=&\sum_{\forall i:\ |B\cap A_i|\le1} (-1)^{|B|}t^{n-|B|}\\
&=&\prod_{i=1}^n(t-|A_i|).\quad\Qed
\eea

To close this section, we point out two situations where our
generalized rank can be described in terms of lengths of chains.

The first of these situations is in a semimodular (hence ranked),
left-modular (hence LL) lattice $L$.  As above, let $\lte$ be the
partial order of $A(L)$ induced by a maximal chain of left-modular
elements and let $\rho$ be the generalized rank function given by the
levels induced by the same chain.  Also, let $\lte^*$ be an arbitrary
linear extension of $\lte$.  We write NBB and NBB${}^*$ to mean NBB
with respect to $\lte$ and $\lte^*$, respectively.  Notice that every
NBB${}^*$ set is also an NBB set.  As noted in our discussion of
Rota's NBC theorem, the size of any NBB${}^*$ base for $x$ is the
ordinary rank of $x$.  By Lemma~\ref{ind}, the size of any NBB base
for $x$ is $\rho(x)$.  So the generalized rank $\rho(x)$ agrees with
the ordinary rank of $x$ provided $x$ has at least one NBB${}^*$ base.
This proviso can be reformulated as $\mu(x)\neq0$,
so the characteristic polynomial is the 
same for both notions of rank.  The proviso
cannot be omitted.  A chain is a distributive lattice (hence  also 
semimodular and left-modular) in which our generalized rank is 1 for all
elements except $\zh$ and therefore differs from the ordinary rank if
the chain has more than two members.

The second situation is described in the following proposition, which
connects $\rho$ to lengths of chains even in some unranked lattices.  
\bpr  			\label{lgch}
Let $L$ be a left-modular lattice which is atomic.  Then for all $x\in
L$ we have
	$$\rho(x)=\mbox{ the length of the longest $\zh$ to $x$ chain.}$$
\epr
\pf\ 
Let $\zh=x_0<x_1<\dots<x_{n-1}<x_n=\oh$ be the left-modular chain used
for the definitions of the sets $A_i$ and thus of $\rho$.  There are,
for any $x\in L$, exactly $\rho(x)$ values of $i$ such that $A_i$
contains an atom $\leq x$.  The elements $x\mt x_i$ for these values
of $i$ (along with $\zh$) constitute a chain of length $\rho(x)$ from
$\zh$ to $x$.  It 
remains to show that no chain from $\zh$ to $x$ is longer, and for
this it suffices to show that, for all $a<b$ in $L$,
$\rho(a)<\rho(b)$.

Let $a<b$ and choose, by atomicity of $L$, an atom $p\leq b$ such that
$p\not\leq a$ and such that $p\in A_j$ for the smallest possible $j$.
Then, by atomicity again, $x_{j-1}\mt b\leq a$.  This and the
left-modularity of $x_{j-1}$ imply 
$$
a=a\jn(x_{j-1}\mt b)=(a\jn x_{j-1})\mt b.
$$
If there were an atom $q\in A_j$ that is $\leq a$, then we would have
$x_{j-1}<q\jn x_{j-1}\leq x_j$ and, by maximality of
the chain of 
$x_i$'s, $q\jn x_{j-1}=x_j$.  But then $(a\jn x_{j-1})\mt b\geq x_j\mt
b\geq p$.  This and $a\not\geq p$ contradict the equation displayed
above.

So there is no such $q$.  But that means that the $A_i$'s counted by
$\rho(b)$ include all those counted by $\rho(a)$ and at least one
more, namely $A_j$.  Therefore, $\rho(a)<\rho(b)$.\hfill\Qed

Notice that the hypotheses of the preceding proposition do not imply
that $L$ is ranked.  A counterexample is given by the six-element
lattice obtained from the eight-element Boolean algebra by removing
two co-atoms.

\section{More examples}					\label{me}

We will now give two examples where our factorization theorem can be
applied but Stanley's cannot because the lattices involved are not
semimodular.  

The first example is the shuffle poset $\cW_{m,1}$
with the maximal chain 
	\beq						\label{Wchain}
	\bx=x_0<x_1<\ldots<x_{m+1}=\by
	\eeq
where $x_1,\ldots,x_m$ are obtained by deleting the letters of $\bx$
in some order.  
(Note that we are using $x_i$ to denote an element of the chain rather
than a letter of $\bx$.)
Greene~\cite{gre:ps} showed that the given chain
satisfies the supersolvability condition
even when extended in $\cW_{m,n}$ by adding the letters of
$\by$ in some order.  So by Proposition~\ref{usmb} $\cW_{m,1}$ is
left-modular. 
We also mentioned in the proof of the same proposition
that this poset satisfies the level condition.
(However the level condition does not hold in $\cW_{m,n}$
for general $n\geq2$, and this is reflected by the fact that 
the corresponding characteristic polynomials usually do not factor
over the integers.)
It is now easy to see that the number of new atoms below $x_i$
in~(\ref{Wchain}) is
	$$|A_i|=\case{1}{if $i\le m$}{m+1}{if $i=m+1$.}$$
{}From Theorem~\ref{LL} we immediately get the following.
\bco
We have
	$$\chi(\cW_{m,1},t)=(t-1)^m(t-m-1).\quad\Qed$$
\eco

Note that $\cW_{m,n}$ is ranked in the ordinary sense and Greene
computed $\chi(\cW_{m,n},t)$ using the usual rank function.  But, 
by Proposition~\ref{lgch}, this rank function coincides with ours.  

For our second example we will use the Tamari 
lattices~\cite{er:hst,ft:tfi,gey:tl,gra:lt,ht:spl}.
Consider all proper parenthesizations $\pi$  of the word
$x_1 x_2\ldots x_{n+1}$.  It is well known that the number of these is
the Catalan number $C_n$.  Partially order
this set by saying that $\si$  covers  $\pi$ whenever
	$$\pi=\ldots ((AB)C)\ldots\qquad\mbox{and}\qquad
	\si=\ldots (A(BC))\ldots$$
for some subwords $A,B,C$.  The corresponding poset turns out to be a
lattice called the {\sl Tamari lattice\/} 
$T_n$. Figure~\ref{T3}(a) gives a picture of $T_3$.

\setlength{\unitlength}{1.3pt}
\bfi
\btab{ll}
\bpi(150,80)(-75,-40)
\Ha
\put(0,35){\mbb{$(x_1(x_2(x_3x_4)))$}}
\Hb
\put(-30,15){\mbr{$(x_1((x_2x_3)x_4))$ }}
\Hc
\put(-30,-15){\mbr{$((x_1(x_2x_3))x_4)$ }}
\Hd
\put(0,-35){\mbt{$(((x_1x_2)x_3)x_4)$}}
\He
\put(30,-15){\mbl{ $((x_1x_2)(x_3x_4))$}}
\Hab
\Hbc
\Hcd
\Hde
\Hae
\epi
&
\bpi(150,80)(-75,-40)
\Ha
\put(0,35){\mbb{$(1,2,3)$}}
\Hb
\put(-30,15){\mbr{$(1,2,2)$ }}
\Hc
\put(-30,-15){\mbr{$(1,2,1)$ }}
\Hd
\put(0,-35){\mbt{$(1,1,1)$}}
\He
\put(30,-15){\mbl{ $(1,1,3)$}}
\Hab
\Hbc
\Hcd
\Hde
\Hae
\epi\\
&\\
(a) Parenthesized version
&
(b) Left bracket version
\etab
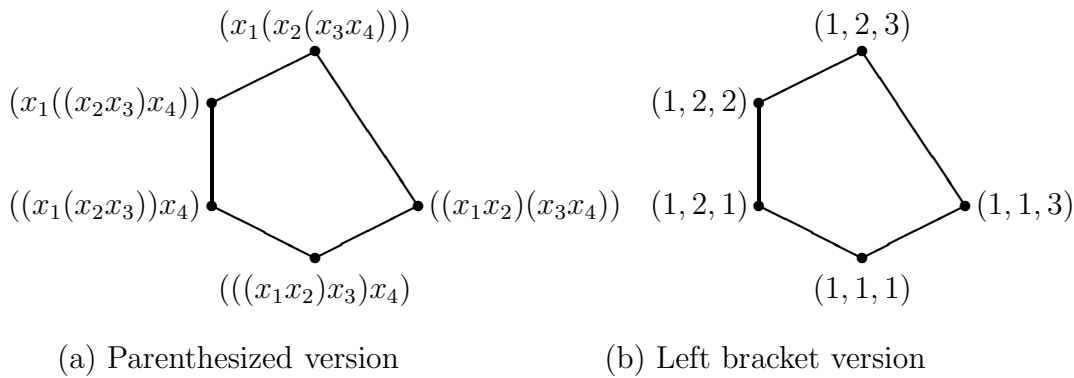
\capt{The Tamari lattice $T_3$}		\label{T3}
\efi

A {\sl left bracket vector\/}, $(v_1,\ldots,v_n)$, is a vector of
positive integers
satisfying
\ben
\item $1\leq v_i\leq i$ for all $i$ and
\item if $S_i=\{v_i,v_i+1,\ldots,i\}$ then for any pair $S_i,S_j$ either
one set contains the other or $S_i\cap S_j=\emp$.
\een
The number of left bracket vectors having
$n$ components is also $C_n$.  In fact given a parenthesized word
$\pi$ we have an associated left bracket vector
$v=(v_1,\ldots,v_n)$ defined as follows.  To calculate $v_i$,
start at $x_i$ in $\pi$ and move left, counting the number of
$x$'s and the number of left parentheses you pass (including $x_i$
itself) until these two
numbers are equal.  Then $v_i=j$ where $x_j$ is the last $x$ 
passed before the numbers balance.  It is not hard to show that this
gives a bijection between parenthesizations and left bracket vectors,
thus inducing a partial order on the latter.  In fact this induced
order is just the component-wise one.  
Figure~\ref{T3}(b) gives the bracket vector version of $T_3$.

Left bracket vectors are also directly related to the trees $G_B$
described in Theorem~\ref{NCNBB}, but with $n+1$ in place of $n$.
Indeed, given a left bracket vector $(v_1,\dots,v_n)$, we obtain such
a tree $G_B$ with vertex set $[n+1]$ by joining $i+1$ to $v_i$ for
$i=1,2,\dots,n$, and all trees as in Theorem~\ref{NCNBB} can be
obtained in this way.

We should note that we have used the left bracket vector rather than
the more traditional right bracket one because when using the former
it is easier to describe the join operation than the meet.  (The
situation is reversed for the right bracket vector.)  This makes it
simpler to work with some of our conditions which only involve joins.
Expressions for the join and meet can be obtained by dualizing results
in~\cite{ht:spl} and~\cite{mar:pie} respectively.
\bpr						\label{jnmt}
Given left bracket vectors $v=(v_1,\ldots,v_n)$ and
$w=(w_1,\ldots,w_n)$ then	
	$$v\jn w=(\max\{v_1,w_1\},\ldots,\max\{v_n,w_n\}).$$
If we let $m=(m_1,\ldots,m_n)$ where $m_i=\min\{v_i,w_i\}$ then 
$v\mt w=(l_1,\ldots,l_n)$ where the $l_i$ are computed recursively by
	$$l_i=\min\{m_i,l_{m_i},l_{m_i+1},\ldots,l_{i-1}\}.\quad\Qed$$
\epr

To show that $T_m$ is LL, consider the chain
	$$\barr{l}
	\De:\ (1,\ldots,1)<(1,2,1,\ldots,1)<(1,2,2,1,\ldots,1)\\
\rp{40}{0}<(1,2,3,1,\ldots,1)<(1,2,3,2,1,\ldots,1)<\ldots<(1,2,3,\ldots,n).
	\earr$$
It is easy to see that this is of maximum length in $T_n$.
The description of join in Proposition~\ref{jnmt} also gives a quick
proof that the level condition holds.  To verify left-modularity will
take more work.

A typical element of $\De$ looks like
	\beq						\label{x}
	x=(1,2,\ldots,j-1,i,1,\ldots,1)
	\eeq
where $i\le j$.  Take $y=(y_1,\ldots,y_n)\le z=(z_1,\ldots,z_n)$ in $T_n$. 
We wish to compute $y\jn(x\mt z)=(c_1,\ldots,c_n)$ so we first
consider $x\mt z$.  Following the 
notation of Proposition~\ref{jnmt} we have
	$$m=(z_1,\ldots,z_{j-1},\min\{z_j,i\},1,\ldots,1).$$
Using the recursive construction on the first $j-1$ components of
$m$ leaves them unchanged since these are the initial components of a
vector in $T_{j-1}$.  Also the last $n-j$ components will still be 1
because $m_i$ is in the minimum taken for $l_i$.  So
$x\mt z$ is the same as $m$ except possibly in the $j$th component
which is
$z_j$ if $z_j\le i$ or $\min\{i,z_i,\ldots,z_{j-1}\}$ if $z_j\ge i$.
Since $y\le z$ we have $y_j\le z_j$ for all $j$ and so
	\beq						\label{y,xz}
	y\jn(x\mt z)=(z_1,\ldots,z_{j-1},c_j,y_{j+1},\ldots,y_j)
	\eeq
where
	\beq						\label{cj}
	c_j=
\case{z_j}{if $z_j\le i$}{\max\{y_j,\min\{i,z_i,\ldots,z_{j-1}\}\}}
{if $z_j\ge i$.}
	\eeq

Now we concentrate on $(y\jn x)\mt z=(d_1,\ldots,d_n)$. Clearly
	$$y\jn x=(1,2,\ldots,j-1,\max\{i,y_j\},y_{j+1},\ldots,y_n)$$
Taking the minimum vector to compute $(y\jn x)\mt z$ we get
	$$m'=(z_1,\ldots,z_{j-1},d_j',y_{j+1},\ldots,y_n)$$
where $d_j'$ is $z_j$ if $z_j\le i$ and $\max\{i,y_j\}$ if $z_j\ge i$.
So in $(y\jn x)\mt z$ we have  $d_i=z_i=c_i$ for $i<j$ by the
same reasoning as before.  We claim that $d_i=y_i$ for $i>j$.  From
the minimization procedure in Proposition~\ref{jnmt} we have 
$d_i\le m_i=y_i$.  But in any lattice $(y\jn x)\mt z\ge y\jn(x\mt z)$
and so 
	\beq						\label{dc}
	d_i\ge c_i \mbox{ for all } 1\le i\le n.  
	\eeq
Since $c_i=y_i$ for $i>j$
this forces the desired equality.

To complete the proof of left-modularity we must show $d_j=c_j$ and
this breaks up into three cases.  If $z_j\le i$ then we have that
$d_j=z_j=c_j$.  If $y_j\le i\le z_j$ then $d_j'=i$ and
$d_j=\min\{i,z_i,\ldots,z_{j-1}\}$.  Comparison of this last equation
with~(\ref{cj}) gives $d_j\le c_j$ and then~(\ref{dc}) results in
$d_j=c_j$.  Finally if $i\le y_j\le z_j$ then $d_j'=y_j$ and
	$$d_j=\min\{y_j,z_{y_j},\ldots,z_{j-1}\}\le y_j=c_j.$$
Using~(\ref{dc}) again finishes this last case. 

It is now an easy matter to compute the characteristic polynomial.
Note that the atoms of $T_n$ are all elements of the form
$(1,\ldots,1,j,1,\ldots,1)$ with the $j\ge2$ in the $j$th position.  If
$x$ is as in~(\ref{x}) then it covers no new atoms if $i<j$ and
exactly one new atom if $i=j$.  Using Theorem~\ref{LL} this translates
as follows.
\bco
We have
	$$\chi(T_n,t)=t^{^{n-1\choose2}}\ (t-1)^{n-1}.\quad\Qed$$
\eco

\section{A further generalization}

In this section, we briefly describe a generalization of
Theorem~\ref{NBB}.  As before, we work with a finite lattice $L$, but
instead of an additional partial ordering of the atoms we use an
arbitrary function $M$ assigning to each $x\in L\setminus\zh$ a
non-empty set $M(x)$ of atoms $\leq x$.  For comparison with
Theorem~\ref{NBB}, the reader should regard the ordering $\lte$ used
there as inducing the function 
$$
M(x)=\{a\in A(L)\mid a \mbox{ is $\lte$-minimal among atoms }\leq x\}.
$$
But what follows applies also to $M$'s that do not arise from 
partial orderings in this way.

For any set $B$ of atoms, let $S(B)$ be the subset obtained by
deleting from $B$ all members of $M(x)$ for all $x\geq\Jn B$.  By the
{\sl core\/} of $B$, we mean the set obtained by starting with $B$ and
repeatedly applying $S$ until the decreasing sequence $B, S(B),
S(S(B)), \dots$ stabilizes; the final $S^n(B)$ is the core of $B$.
(This can also be described as the largest subset of $B$ unchanged by
$S$.)  We call $B$ {\sl coreless\/} if its core is empty. 
\bth
Let $L$ be any finite lattice and let $M$ be any function assigning to
every $x\in L\setminus\zh$ a nonempty set $M(x)$ of atoms that are $\leq
x$.  Then for all $x\in L$ we have
$$
\mu(x)=\sum_B (-1)^{|B|}
$$
where the sum is over all coreless $B\subseteq A(L)$ whose join is
$x$. 
\eth

We omit the proof since it is essentially the same as that of
Theorem~\ref{NBB}.  The only difference is that, instead of choosing a
$\lte$-minimal atom $a_0\leq x$, we choose an $a_0\in M(x)$.

Notice that, if $M$ happens to be obtained from a partial order $\lte$
as described above, then non-empty cores with respect to $M$ are the
same as BB sets with respect to $\lte$, and therefore coreless sets
with respect to $M$ are the same as NBB sets with respect to $\lte$.

Notice also that the theorem of this section is ``stronger,'' in the
sense of having fewer summands in the formula for the M\"obius
function, when $M(x)$ is smaller, for this will make $S(B)$ larger and
therefore make $B$ less likely to be coreless.  From this point of
view, one should always take $M(x)$ to be a singleton.  Of course
considerations of naturality or clarity may make other choices of $M$
preferable (just as they may make nonlinear orderings $\lte$ preferable
to the more efficient linear orderings in applications of
Theorem~\ref{NBB}).  

\section{Comments and Questions}

\hspace{\parindent} (1)  It would be interesting to find other
applications of our two 
main theorems.    The higher Stasheff-Tamari posets as recently
defined in~\cite{er:hst} are obvious candidates.  However, one would
first have to prove that they are lattices.

\medskip

(2)  In a previous paper~\cite{sag:grn} one of us proved a
somewhat weaker generalization of Rota's NBC Theorem which,
nonetheless, has some interesting connections to our work and others.
To state the result we must first recall some definitions
from~\cite{sag:grn}.  Call  $B\sbe A(L)$ {\sl independent\/}
if $\Jn B'<\Jn B$ for every proper subset $B'\sbs B$.  So $C$ is
dependent if $\Jn C' = \Jn C$ for some $C'\sbs C$.  Note that these
generalize the corresponding notions for a geometric lattice.  The
definitions of base, circuit, and broken circuit remain as before.  The
main result of~\cite{sag:grn} can now be stated.
\bth						\label{oldNBC}
Let $L$ be a finite lattice. Let $\lte$ be any total ordering of
$A(L)$ such that for all ciruits $C$ we have
	\beq					\label{C'}
	\Jn C =\Jn (C\setm\min C).
	\eeq
Then for all $x\in L$ we have
	$$\mu(x)=\sum_B (-1)^{|B|}$$
where the sum is over all NBC bases of $x$.\hfill\Qed
\eth

We should note that not every lattice has a total ordering of the
atoms satisfying~(\ref{C'}) and so this result is not as strong as
Theorem~\ref{NBB}.  On the other hand, there is an interesting
relationship between a dual of this restriction and the level condition from
Section~\ref{sl}. 
\bpr
Let $L$ be a finite lattice. Let $\lte$ be an induced ordering of
$A(L)$ such that for all ciruits $C$ having a unique maximal element 
we have
	$$\Jn C =\Jn (C\setm\max C).$$
Then $L$ satisfies the level hypothesis.
\epr
\pf\  Suppose, toward a contradiction, that we have 
$b_0\lt b_1\lt\ldots\lt b_k$ and $b_0\le\Jn_{i\geq1} b_i$.  Then
$D=\{b_0,b_1,\ldots,b_k\}$ is dependent and so contains a circuit $C$.
Also $C$ must have a unique maximal element $b_j$ because $D$ intersects
each level at most once.  So we have $\Jn C = \Jn (C\setm b_j)$, or
equivalently $b_j\le\Jn (C\setm b_j)$.  But this cannot happen since
for some $x$ of the inducing chain we have $b_i\le x$ for $i<j$ but
$b_j\not\le x$.\hfill\Qed

In his work~\cite{ath:tph} on Tutte polynomials for hypermatroids, Christos
Athanasiadis defines a {\sl generalized lattice\/} to be a triple
$\cL=(L,E,f)$ where $L$ is a finite lattice, $E$ is a set, and $f:E\ra L$
is some map.  He then proves a theorem about the M\"obius function of
$L$ in the case that the hypermatroid associated with $\cL$ satifies
condition~(\ref{C'}).  In the case that $E=A(L)$ and $f$ is inclusion
his result becomes a special case of Theorem~\ref{oldNBC}.

\medskip

(3)  There are certain to be topological ramifications of our work.
In fact Yoav Segev~\cite{seg:scn} has already proved that the order
complex of any lattice is homotopy equivalent to the simplicial
complex of all NBB sets $B$ with $\Jn B<\oh$.  This can be used to
demonstrate a recent result of Linusson~\cite{lin:cli} that the order
complex of an interval in the partition lattice under dominance 
is homotopy equivalent to a sphere or contractible.

An important use of NBC sets is to give a basis for the Orlik-Solomon
algebra $\bA(L)$ of a geometric lattice $L$~\cite{os:ctc}.  If $L$ is
the intersection lattice of a complex hyperplane arrangement, then $\bA(L)$ is
isomorphic to 
the cohomology algebra of the complement of the arrangement.
Recently there has been a lot of interest in subspace
arrangements~\cite{bjo:sa}.  The corresponding intersection lattices
are no longer geometric, but perhaps NBB sets can be used to 
give information about the associated cohomology algebra in this case.
Recently, De Concini and Procesi~\cite{dp:wms} used algebraic geometric
techniques 
to show that this algebra is indeed determined solely by the lattice
and dimension information.  This provides some hope that a
combinatorial description is also possible.

\begin{\bib}{99}

\bibitem{ath:tph} C. A. Athanasiadis, The Tutte polynomial of a
hypermatroid, preprint, 1995.

\bibitem{bjo:sa}  A. Bj\"orner, Subspace arrangements, in ``Proc. 1st
European Congress Math. (Paris 1992),'' A. Joseph and R. Rentschler
eds., Birkh\"auser, Boston, MA, to appear.

\bibitem{bl:ldt}  A. Bj\"orner and L. Lov\'asz,
Linear decision trees, subspace arrangements and M\"obius functions,
{\sl \jams} {\bf 7} (1994), 677--706.

\bibitem{bly:ldt}  A. Bj\"orner, L. Lov\'asz and A. Yao,
Linear decision trees: volume estimates and topological bounds, in
``Proc. 24th ACM Symp. on Theory of Computing,''  ACM Press, New York,
NY, 1992, 170--177.

\bibitem{bog:mfd} K. Bogart, The M\"obius function of the domination
lattice, unpublished manuscript, 1972.

\bibitem{bry:lip}  T. Brylawski, The lattice of integer partitions,
{\sl \dm} {\bf 6} (1973), 201--219.

\bibitem{dp:wms} C. De Concini and C. Procesi, Wonderful models of
subspace arrangements, preprint, 1995.

\bibitem{er:hst} P. H. Edelman and V. Reiner, The higher Stasheff-Tamari
posets, preprint.

\bibitem{ft:tfi} H. Friedman and D. Tamari, Probl\`emes
d'associativit\'e:  Une treillis finie induite par une loi
demi-associative, {\it J. Combin. Theory} {\bf 2} (1967), 215--242.

\bibitem{gey:tl} W. Geyer, On Tamari lattices, {\sl \dm} {\bf 133}
(1994), 99--122.

\bibitem{gra:lt} G. Gr\"atzer, ``Lattice Theory,'' Freeman and Co.,
San Francisco, CA, 1971, pp. 17--18, problems 26-36.

\bibitem{gre:clm} C. Greene, A class of lattices with M\"obius
function $\pm1$, {\sl \ejc} {\bf 9} (1988), 225--240.

\bibitem{gre:ps} C. Greene, Posets of Shuffles, {\sl \jcta} {\bf 47}
(1988), 191--206.

\bibitem{ht:spl} S. Huang and D. Tamari, Problems of associativity:  A
simple proof for the lattice property of systems ordered by a
semi-associative law, {\it \jcta} {\bf 13} (1972), 7--13.

\bibitem{kre:pnc}  G. Kreweras, Sur les partitions non-crois\'ees d'un
cycle, {\sl \dm} {\bf 1} (1972), 333--350.

\bibitem{lin:cli}  S. Linusson, A class of lattices whose intervals
are spherical or contractible, preprint.

\bibitem{mar:pie} G. Markowsky, Primes, irreducibles and extremal
lattices, {\sl Order} {\bf 9} (1992), 7--13.

\bibitem{os:ctc} P. Orlik and L. Solomon, Combinatorics and topology
of complements of hyperplanes. {\sl Invent. Math.} {\bf 56} (1980),
167--189.

\bibitem{rei:npc}  V. Reiner, Non-crossing partitions for classical
reflection groups, preprint.

\bibitem{rot:tmf} G.-C. Rota, On the foundations of combinatorial
theory I. Theory of M\"obius functions, {\sl Z. 
Wahrscheinlichkeitstheorie} {\bf 2} (1964), 340--368.

\bibitem{sag:grn} B. E. Sagan, A generalization of Rota's NBC Theorem,
{\sl Adv. in Math.}, {\bf 111} (1995), 195--207.

\bibitem{seg:scn}  Y. Segev, The simplicial complex of all NBB
nonspanning subsets of the set of atoms of a prelattice $L$ is
homotopic to the order complex of $L$, preprint, 1995.

\bibitem{sta:ssl} R. P. Stanley, Supersolvable lattices, {\sl
Alg. Univ.}, {\bf 2} (1972), 197--217.

\bibitem{sta:ec1} R. P. Stanley, ``Enumerative Combinatorics,
Volume 1,''  Wadsworth \& Brooks/Cole, Pacific Grove, CA, 1986.

\bibitem{whi:lem}  H. Whitney, A logical expansion in mathematics,
{\sl Bull. Amer. Math. Soc.} {\bf 38} (1932), 572--579.

\bibitem{zas:grs} T. Zaslavsky, The geometry of root systems
and signed graphs, {\it Amer. Math. Monthly} {\bf 88} (1981), 88--105.

\end{\bib}

\end{document}